\input amstex
\documentstyle{amsppt}
%
%
\nopagenumbers
\accentedsymbol\tx{\tilde x}
\def\Img{\operatorname{Im}}
\def\Ker{\operatorname{Ker}}
\def\id{\operatorname{id}}
\def\tr{\operatorname{tr}}
\def\negskp{\hskip -2pt}
\def\compos{\,\raise 1pt\hbox{$\sssize\circ$} \,}
\pagewidth{360pt}
\pageheight{606pt}
\leftheadtext{Ruslan A. Sharipov}
\rightheadtext{On the concept of normal shift \dots}
\topmatter
\title On the concept of normal shift in non-metric geometry.
\endtitle
\author
R.~A.~Sharipov
\endauthor
\abstract
Theory of Newtonian dynamical systems admitting normal shift of
hypersurfaces was first developed for the case of Riemannian
manifolds. Recently it was generalized for manifolds geometric
equipment of which is given by some regular Lagrangian or,
equivalently, by some regular Hamiltonian dynamical system. In
present paper we consider further generalization of this theory
for the case, when geometry of manifold is given by generalized
Legendre transformation.
\endabstract
\address Rabochaya street 5, 450003, Ufa, Russia
\endaddress
\email \vtop to 20pt{\hsize=280pt\noindent
R\_\hskip 1pt Sharipov\@ic.bashedu.ru\newline
r-sharipov\@mail.ru\vss}
\endemail
\urladdr
http:/\negskp/www.geocities.com/r-sharipov
\endurladdr
\endtopmatter
\loadbold
\TagsOnRight
\document
\head
1. What is normal shift\,? Brief historical overview.
\endhead
\parshape 4 0pt 360pt 0pt 360pt 0pt 360pt 180pt 180pt
    Phenomenon of normal shift is very simple by its nature. Let's
consider it in three-dimensional Euclidean space $\Bbb R^3$. Suppose
that $\sigma$ is some smooth orientable surface in $\Bbb R^3$. At each
point $p\in \sigma$ one can draw unit normal vector $\bold n$ such that
$\bold n=\bold n(p)$ would be a smooth vector-valued function on $\sigma$.
Let's move each point $p$ of $\sigma$ in the direction of vector $\bold
n(p)$ to the distance $t$ which is the same for all points $p\in \sigma$.
\vadjust{\vskip -46pt\hbox to 0pt{\kern 5pt\hbox{\special{em:graph
pst-17a.gif}}\hss}\vskip 46pt}Then moved points $p_t$ would form
another surface $\sigma_t$ as shown on Fig\. 1.1. Changing parameter $t$
we would obtain one-parametric family of surfaces. This construction
is known as Bonnet transformation.\par
\parshape 5 180pt 180pt 180pt 180pt 180pt 180pt 180pt 180pt
0pt 360pt
     In Bonnet construction initial surface $\sigma$ is transformed by
moving each point of $\sigma$. Trajectories of motion in this case are 
straight
lines directed along normal vectors and points of $\sigma$ move along them
with a constant speed $|\bold v|=1$. Therefore parameter $t$, which is
the distance of displacement, can also be interpreted as time variable.
Bonnet noted that all surfaces $\sigma_t$ in his construction are
perpendicular to the trajectories of moving points. For this reason his
construction is also known as {\bf normal displacement} or {\bf normal
shift}.\par
    Basic observation by Bonnet, i\.\,e\. orthogonality of surfaces
$\sigma_t$ and shift trajectories, gave an impetus for generalization of
his construction. This was done by me and my student A.~Yu\.~Boldin
in preprint \cite{1} (see also \cite{2}, \cite{3}). We replaced straight
lines in Bonnet construction by trajectories of Newtonian dynamical system:
$$
\xalignat 2
&\hskip -2em
\dot\bold r=\bold v,
&&\dot\bold v=\bold F(\bold r,\bold v).
\tag1.1
\endxalignat
$$
This dynamical system describes the motion of a particle of unit mass
according to Newton's second law. Vector $\bold F$ in right hand side
of \thetag{1.1} is a vector of force acting on this particle. In order
to initialize a shift of surface $\sigma$ we applied initial data
$$
\xalignat 2
&\hskip -2em
\bold r\,\hbox{\vrule height 8pt depth 8pt width 0.5pt}_{\,t=0}
=\bold r(p),
&&\bold v\,\hbox{\vrule height 8pt depth 8pt width 0.5pt}_{\,t=0}
=\nu(p)\cdot\bold n(p)
\tag1.2
\endxalignat
$$
to differential equations \thetag{1.1}. For classical Bonnet
transformation $|\bold v|=1$. In our construction modulus of velocity
vector is not constant. It's initial value on $\sigma$ is determined by
some scalar function $\nu=\nu(p)$.\par
    Initial data \thetag{1.2} with parameter $p\in \sigma$ determine a
family of trajectories of Newtonian dynamical system \thetag{1.1}. 
Points of $\sigma$ moving along these trajectories form one-parametric
family of surfaces $\sigma_t$. Thus we have generalization of Bonnet
construction. This is the shift of $\sigma$ along trajectories of Newtonian
dynamical system \thetag{1.1}. It's clear that this is normal shift
for initial instant of time $t=0$. However, in general, orthogonality
of $\sigma_t$ and shift trajectories gets broken at any other instant of time
$t\neq 0$. Only for special Newtonian dynamical systems, i\.\,e\. for
special force fields $\bold F=\bold F(\bold r,\bold v)$, one can keep
orthogonality for $t\neq 0$ at the expense of proper choice of function
$\nu=\nu(p)$ in \thetag{1.2}.
\definition{Definition 1.1} Shift of surface $\sigma\subset\Bbb R^3$ along
trajectories of dynamical system \thetag{1.1} determined by initial
data \thetag{1.2} is called {\bf normal shift} if all displaced surfaces
$\sigma_t$ are perpendicular to shift trajectories.
\enddefinition
\definition{Definition 1.2} We say that Newtonian dynamical system
\thetag{1.1} satisfies {\bf normality condition} if for sufficiently small
part of any surface $\sigma\subset\Bbb R^3$ there is a function $\nu=\nu(p)$ 
on
it determining normal shift of this part of $\sigma$.
\enddefinition
    Words ``sufficiently small part'' of $\sigma$ in definition~1.2 mean
that for any fixed point $p_0\in \sigma$ there is some sufficiently small
open neighborhood of $p_0$ in $\sigma$, where proper function $\nu=\nu(p)$
does exist. Let $\nu_0\neq 0$ be some arbitrary nonzero constant.
Then we can normalize function $\nu(p)$ at the fixed point $p_0$:
$$
\hskip -2em
\nu(p_0)=\nu_0.
\tag1.3
$$
\definition{Definition 1.3} Newtonian dynamical system \thetag{1.1}
satisfies {\bf strong normality condition} if for any surface $\sigma\subset
\Bbb R^3$, for any fixed point $p_0\in \sigma$, and for any constant $\nu_0
\neq 0$ there is some open neighborhood of $p_0$ and some smooth
function $\nu(p)$ normalized by the condition \thetag{1.3} in this
neighborhood such that it determines normal shift in the sense of
definition~1.1.
\enddefinition
    Newtonian dynamical systems satisfying strong normality condition
form special subclass which appears to be interesting object for study.
We call them {\bf systems admitting normal shift} of hypersurfaces. In
simpler words, these are systems capable to implement normal displacement
of any hypersurface $\sigma$ in $\Bbb R^3$ with any predefined value $\nu_0$
of initial velocity.\par
    It's obvious that definitions~1.1, 1.2, and 1.3 can be formulated for
higher dimensional Euclidean spaces and for Riemannian manifolds as well.
One should only replace surfaces by hypersurfaces and replace $\Bbb R^3$
by $\Bbb R^n$ or by arbitrary smooth manifold $M$ with Riemannian metric
$\bold g$. This was done in papers \cite{4--6}. Instead of \thetag{1.1}
in Riemannian manifolds we write differential equations
$$
\xalignat 2
&\hskip -2em
\dot x^i=v^i,
&&\nabla_{\!t}v^i=F^i(x^1,\ldots,x^n,v^1,\ldots,v^n),
\tag1.4
\endxalignat
$$
where $i=1,\,\ldots,\,n$. Here $x^1,\,\ldots,\,x^n$ are coordinates
of a point $p$ in some local chart of Riemannian manifold $M$, while
$v^1,\,\ldots,\,v^n$ are components of velocity vector $\bold v\in
T_p(M)$. Vector $\bold F\in T_p(M)$ with components $F^1,\,\ldots,\,
F^n$ determines force field of Newtonian dynamical system \thetag{1.4}.
\par
    In papers \cite{5} and \cite{6} we have shown that strong normality
condition applied to Newtonian dynamical system \thetag{1.4} leads to a
system of partial differential equations for components of force vector
$\bold F=\bold F(p,\bold v)$. This system subdivides into two parts:
\roster
\item"--" {\bf weak normality equations} written for $n\geqslant 2$;
\item"--" {\bf additional normality equations} written for $n\geqslant 3$.
\endroster
Note that $n=2$ is lower limit for the dimension of manifold $M$. Indeed,
for $n=1$, hypersurfaces are points, concept of normal shift in this case
has no meaning. Note also that additional normality equations are written
for $n\geqslant 3$. Therefore $n=2$ is exceptional dimension. Theory of
dynamical systems admitting normal shift in two-dimensional case $n=2$ is
rather different from that of multidimensional case $n\geqslant 3$. This
fact is reflected in theses \cite{7} and \cite{8}.\par
    In two-dimensional case strong normality condition is equivalent to
weak normality equations for force field $\bold F$. Weak normality
equations in this case can be reduced to one nonlinear partial differential
equation for one scalar function of four variables. This equation cannot be
solved explicitly in general. However, one can construct some special
explicit solutions of it (see paper \cite{9} and thesis \cite{8}, where
numerous examples are given).\par
    In multidimensional case $n\geqslant 3$ strong normality condition is
equivalent to complete system of weak normality equations and additional
normality equations for components force field $\bold F$. As appeared, in
this case one can find an explicit formula for general solution of this
rather huge system of PDE's (see papers \cite{10}, \cite{11}, and
Chapter~\uppercase\expandafter{\romannumeral 7} of thesis \cite{7}).
Here is this formula:
$$
\hskip -2em
F_i=\frac{h(W)}{W_v}\cdot\frac{\,v_i}{|\bold v|}-\sum^n_{k=1}
\frac{\nabla_kW}{W_v}\cdot\frac{2\,v^k\,v_i-|\bold v|^2\,
\delta^k_i}{|\bold v|}.
\tag1.5
$$
Formula \thetag{1.5} contain two arbitrary functions $W$ and $h$,
where $h=h(w)$ is an arbitrary smooth function of one variable,
while $W=W(x^1,\ldots,x^n,v)$ is an arbitrary smooth function of
$n+1$ variables with nonzero derivative
$$
W_v=\frac{\partial W}{\partial v}\neq 0.
$$
By $\nabla_kW$ in formula \thetag{1.5} we denote partial derivatives
$$
\nabla_kW=\frac{\partial W}{\partial x^k},
$$
while $v_i$ and $v^k$ are covariant and contravariant components of
velocity vector $\bold v$:
$$
v_i=\sum^n_{k=1}g_{ik}\,v^k.
$$
When substituted into \thetag{1.5}, last $(n+1)$-th argument $v$
of $W_v$ and $\nabla_kW$ is replaced by modulus of velocity vector:
$v=|\bold v|$.\par
    Division of complete system of normality equations into two parts
is not artificial by its nature. As appeared, weak normality equations
have their own geometrical interpretation. In paper \cite{12} normal
blow-up of points in Riemannian manifolds was considered. There it was
found that force field of Newtonian dynamical systems admitting normal
blow-up of points should satisfy weak normality equations only. This
means that in multidimensional case $n\geqslant 3$ they could form
larger class of dynamical systems than those admitting normal shift of
hypersurfaces. Examples given in paper \cite{13} show that they actually
do form larger class.\par
    Papers \cite{14} and \cite{15} are devoted to the study of global
geometric structures associated with Newtonian dynamical systems
admitting normal shift of hypersurfaces in multidimensional case
$n\geqslant 3$. This study is based mainly on explicit formula
\thetag{1.5}. In addition, one should note that theory of dynamical
systems admitting normal shift of hypersurfaces was generalized to
the case of Finslerian manifolds
(see Chapter~\uppercase\expandafter{\romannumeral 8} of thesis \cite{7}).
Weak and additional normality equations were derived. However, explicit
formula like \thetag{1.5} for this case is not yet obtained.\par
    During one year after the conference dedicated to Centenary Anniversary
of I.~G.~Petrovsky, May 2001, I was looking for applications of the theory
constructed. Results are represented by papers \cite{16} and \cite{17}.
The idea is very simple. It is known that in the limit of short waves all
wave propagation phenomena can be described in terms of rays and beams
(geometrical optics is an example). In this limit amplitude of scalar wave
is described by asymptotic expansion
$$
\hskip -2em
u=\sum^\infty_{\alpha=0}\frac{\varphi_{\sssize(\ssize\alpha\sssize)}}
{(i\,\lambda)^\alpha}\cdot e^{i\lambda S},\quad\lambda\to\infty,
\tag1.6
$$
known as Debye's ansatz (see \cite{18}). Here in \thetag{1.6} function
$S=S(x^1,\ldots,x^n)$ is a phase of propagating wave. This function
satisfies Hamilton-Jacobi equation
$$
\hskip -2em
H(x^1,\ldots,x^n,\nabla_1S,\ldots,\nabla_nS)=0.
\tag1.7
$$
Here $\nabla_1S,\,\ldots,\,\nabla_nS$ are components of momentum covector
$\bold p=\nabla S$:
$$
\hskip -2em
p_i=\nabla_iS=\frac{\partial S}{\partial x^i}.
\tag1.8
$$
Hamilton function $H=H(x^1,\ldots,x^n,p_1,\ldots,p_n)$ in \thetag{1.7}
is determined by wave operator describing physical properties of medium,
where wave propagation process occur. Usually it is polynomial with
respect to components of momentum covector \thetag{1.8}. But in more
complicated cases it may be \pagebreak non-polynomial as well.
\par
   Note that Hamilton-Jacobi equation \thetag{1.7} is a first order
PDE with respect to phase function $S=S(x^1,\ldots,x^n)$. Its solution
is written in terms of characteristic lines. In present case they
are given by Hamilton equations
$$
\xalignat 2
&\hskip -2em
\dot x^i=\frac{\partial H}{\partial p_i},
&&\dot p_i=-\frac{\partial H}{\partial x^i}.
\tag1.9
\endxalignat
$$
Physically, characteristic lines or trajectories of Hamiltonian dynamical
system \thetag{1.9} are interpreted as rays or beams in wave propagation
process. {\bf Wave fronts} in this process are hypersurfaces, where phase
is constant. In other words, they are level hypersurfaces of phase function
$S=S(x^1,\ldots,x^n)$. These hypersurfaces represent the position of wave
at various time instants, so we can say that wave propagates moving along
trajectories of Hamiltonian dynamical system \thetag{1.9}. However, time
variable $t$, with respect to which ordinary differential equations
\thetag{1.9} are written, is not an actual time in wave propagation process.
Actual time is proportional to the value of phase itself, therefore we
can take $t=S$. As shown in \cite{17}, passing to this new time variable,
we obtain the following modified Hamilton equations:
$$
\xalignat 2
&\hskip -2em
\dot x^i=\frac{1}{\Omega}\cdot\frac{\partial H}{\partial p_i},
&&\dot p_i=-\frac{1}{\Omega}\cdot\frac{\partial H}{\partial x^i}.
\tag1.10
\endxalignat
$$
Denominator $\Omega$, which is often interpreted as kinetic energy,
is given by formula
$$
\hskip -2em
\Omega=\sum^n_{k=1}\frac{\partial H}{\partial p_k}\,p_k.
\tag1.11
$$\par
    Now suppose that we have Riemannian manifold $M$. Then we can
define modulus of momentum covector $p=|\bold p|$ and can consider
Hamilton function of special form $H=W(x^1,\ldots,x^n,|\bold p|)$.
Modified Hamiltonian dynamical system \thetag{1.10} with such special
Hamilton function appears to be equivalent to Newtonian dynamical
system \thetag{1.4} with force field given by the following formula:
$$
\hskip -2em
F_i=-\sum^n_{k=1}\frac{\nabla_kW}{W_v}\cdot\frac{2\,v^k\,v_i
-|\bold v|^2\,\delta^k_i}{|\bold v|}.
\tag1.12
$$
This is the main result of paper \cite{17}. Comparing \thetag{1.5}
and \thetag{1.12} we conclude that part of Newtonian dynamical systems
admitting normal shift of hypersurfaces in Riemannian geometry can be
interpreted as the equations of wave front dynamics in physics. And
conversely, for some wave propagation phenomena we observe something
like {\bf conservation law}: wave front dynamics {\bf preserves
orthogonality} of wave fronts and rays. However, unlike mathematical
theorems, true laws of nature cannot be sharply specific, i\.\,e\.
applicable to some objects and not applicable to others. One should
expect that the above {\bf orthogonality law} can be generalized for
all wave propagation phenomena.\par
    Now suppose that modified Hamilton equations \thetag{1.10} are
written for arbitrary smooth manifold $M$. In the absence of metric
one should revise the concept of orthogonality itself. This was done
in paper \cite{19}. Note that modified Hamilton equations \thetag{1.10}
define dynamics in cotangent bundle $T^*\!M$. Therefore at each point
of trajectory $p=p(t)$ in $M$ we have momentum covector $\bold p=
\bold p(t)$. Therefore if trajectory $p=p(t)$ crosses some hypersurface
$\sigma$ and if $\boldsymbol\tau$ is a vector tangent to $\sigma$, then
we can consider scalar product of $\boldsymbol\tau$ and momentum covector
$\bold p$:
$$
\hskip -2em
\varphi=\left<\bold p\,|\,\boldsymbol\tau\right>
=\sum^n_{k=1}\tau^k\,p_k.
\tag1.13
$$
We see that scalar product \thetag{1.13} does not require any metric.
If $\varphi=0$ for all $\boldsymbol\tau\in T_p(\sigma)$, then we say
that trajectory $p=p(t)$ is perpendicular to hypersurface $\sigma$.
\par
    Let $\sigma$ be some arbitrary orientable hypersurface in $M$.
In the absence of metric we cannot define normal vector for $\sigma$.
However, at each point $p$ of $\sigma$ there is normal covector
$\bold n=\bold n(p)$ perpendicular to $\sigma$ in the sense of scalar
product \thetag{1.13}. It is defined uniquely up to a scalar factor:
$\bold n\to\alpha\cdot\bold n$. We can choose $\bold n=\bold n(p)$
to be smooth covector-valued function on $\sigma$. Then we can define
the following initial data for ordinary differential equations
\thetag{1.10}:
$$
\xalignat 2
&\hskip -2em
x^i\,\hbox{\vrule height 8pt depth 8pt width 0.5pt}_{\,t=0}
=x^i(p),
&&p_i\,\hbox{\vrule height 8pt depth 8pt width 0.5pt}_{\,t=0}
=\nu(p)\cdot n_i(p)
\tag1.14
\endxalignat
$$
This defines motion of the points of $\sigma$ along trajectories of
dynamical system \thetag{1.10}.
\definition{Definition 1.4} Shift of hypersurface $\sigma\subset M$
along trajectories of a dynamical system in cotangent bundle $T^*\!M$
determined by initial data \thetag{1.14} is called {\bf normal shift}
if all displaced hypersurfaces $\sigma_t$ are perpendicular to shift
trajectories in the sense of scalar product determined by formula
\thetag{1.13}.
\enddefinition
\definition{Definition 1.5} Dynamical system in cotangent bundle
$T^*\!M$ satisfies {\bf strong normality condition} if for any
hypersurface $\sigma\subset M$, for any fixed point $p_0\in \sigma$,
and for any constant $\nu_0\neq 0$ there is some open neighborhood
of $p_0$ and some smooth function $\nu(p)$ normalized by the condition
\thetag{1.3} in this neighborhood such that it determines normal shift
of $\sigma$ in the sense of definition~1.4.
\enddefinition
    In paper \cite{19} it was shown that any modified Hamiltonian
dynamical system \thetag{1.10} with arbitrary Hamilton function $H$
(provided $\Omega\neq 0$) satisfies strong normality condition.
This means that dynamical system \thetag{1.10} form background for
studying strong normality condition. In Riemannian geometry such
background is given by geodesic flows (these are Newtonian dynamical
systems \thetag{1.4} with identically zero force field $\bold F=0$).
In paper \cite{19} I considered the following class of dynamical
systems in cotangent bundle $T^*\!M$:
$$
\xalignat 2
&\hskip -2em
\dot x^i=\frac{1}{\Omega}\cdot\frac{\partial H}{\partial p_i},
&&\dot p_i=-\frac{1}{\Omega}\cdot\frac{\partial H}{\partial x^i}+Q_i.
\tag1.15
\endxalignat
$$
Under sufficiently non-restrictive conditions for Hamilton function
$H$ (see \cite{19}) differential equations \thetag{1.15} can be
transformed to the form similar to \thetag{1.4}:
$$
\xalignat 2
&\hskip -2em
\dot x^i=v^i,
&&\dot v^i=\Phi^i(x^1,\ldots,x^n,v^1,\ldots,v^n).
\tag1.16
\endxalignat
$$
Therefore we say that \thetag{1.15} is {\bf relative form} of Newtonian
dynamical system \thetag{1.16} as related to Hamiltonian system
\thetag{1.10}. Covector $\bold Q$ with components $Q_1,\,\ldots,\,Q_n$
in \thetag{1.15} plays the same role as force vector $\bold F$ in
\thetag{1.4}.\par
   Note that Newtonian dynamical system in relative form \thetag{1.15}
is a dynamical system in cotangent bundle $T^*\!M$. Therefore
definitions~1.4 and 1.5 can be applied to it. This was actually done
in paper \cite{20}. In that paper it was shown that strong normality
condition for Newtonian dynamical system \thetag{1.15} is equivalent
to a system of partial differential equations for components of
covector $\bold Q=\bold Q(p,\bold p)$. This system of partial
differential equations subdivides into two parts:
\roster
\item"--" {\bf weak normality equations} written for $n\geqslant 2$;
\item"--" {\bf additional normality equations} written for $n\geqslant 3$.
\endroster
Studying these newly derived normality equations is rather interesting
problem. However, below we go further and we consider more general
situation, when Hamiltonian and/or Lagrangian dynamical system in $M$
is not given.
\head
2. The idea of further generalization.
\endhead
    In the absence of special geometric structures (like Riemannian
metric or Hamilton function) the only way of defining Newtonian
dynamics in $M$ is given by the equations \thetag{1.16}. Can we use
scalar product \thetag{1.13} for to define normal shift in this case.
The answer is {\bf yes}, provided we have some way to determine
momentum co\-vector $\bold p$\,! For instance, momentum covector
$\bold p$ can be given explicitly as a function of dynamic variables
$p$ and $\bold v$. In local chart this looks like
$$
\hskip -2em
\cases p_1=p_1(x^1,\ldots,x^n,v^1,\ldots,v^n),\\
.\ .\ .\ .\ .\ .\ .\ .\ .\ .\ .\ .\ .\ .\ .\ .\
.\ .\ .\ .\ .\ .\ .\ \\
p_n=p_n(x^1,\ldots,x^n,v^1,\ldots,v^n).\\
\endcases
\tag2.1
$$
If coordinates $x^1,\,\ldots,\,x^n$ of the point $p\in M$ are fixed, then
functions \thetag{2.1} express $n$ variables $p_1,\,\ldots,\,p_n$ through
another set of $n$ variables $v^1,\,\ldots,\,v^n$. Therefore one can treat
\thetag{2.1} as coordinate representation of a map
$$
\hskip -2em
\lambda\!:TM\to T^*\!M.
\tag2.2
$$
This map is an analog of well known Legendre transformation (see
\cite{21}), for which functions \thetag{2.1} are determined by
Lagrange function $L=L(p,\bold v)$. Below we shall call it {\bf
generalized Legendre transformation} and for the sake of convenience
we shall assume this map \thetag{2.2} to be diffeomorphism.\par
    Using generalized Legendre transformation, we can transform
Newtonian dynamical system \thetag{1.16} to new dynamic variables
$p$ and $\bold p$ forming point $q=(p,\bold p)$ of cotangent bundle
$TM$. As a result we get the following differential equations:
$$
\hskip -2em
\left\{\aligned
&\dot x^i=V^i(x^1,\ldots,x^n,p_1,\ldots,p_n),\\
&\dot p_i=\Theta_i(x^1,\ldots,x^n,p_1,\ldots,p_n).
\endaligned\right.
\tag2.3
$$
Here $V^1,\,\ldots,\,V^n$ are functions that implement inversion of
the map \thetag{2.2}:
$$
\hskip -2em
\cases v^1=V^1(x^1,\ldots,x^n,p_1,\ldots,p_n),\\
.\ .\ .\ .\ .\ .\ .\ .\ .\ .\ .\ .\ .\ .\ .\ .\
.\ .\ .\ .\ .\ .\ .\ \\
v^n=V^n(x^1,\ldots,x^n,p_1,\ldots,p_n).\\
\endcases
\tag2.4
$$
They define kinematic part of the equations of Newtonian dynamics
\thetag{2.3}. Their role is similar to the role of Hamilton function
in \thetag{1.15}. Functions $\Theta_1,\,\ldots,\Theta_n$ define
dynamical part of the equations \thetag{2.3}, they play the same role
as functions $Q_1,\,\ldots,\,Q_n$ in \thetag{1.15} and functions
$\Phi^1,\,\ldots,\Phi^n$ in \thetag{1.16}.\par
    Note that Newtonian dynamical system written as \thetag{2.3}
is a dynamical system in cotangent bundle $T^*\!M$. Definitions~1.4
and 1.5 are applicable to it. Therefore we can start new theory of
dynamical systems admitting normal shift of hypersurfaces for
manifolds equipped only with generalized Legendre map. Our main
goal here is to derive weak and additional normality equations for
this case.
\head
3. Newtonian dynamics in manifolds and associated geometric structures.
\endhead
    Let's consider Newtonian dynamical system in some smooth manifold $M$.
In the absence of metric we cannot use covariant derivative $\nabla_{\!t}$
in \thetag{1.4}. Therefore we write it as \thetag{1.16}. This is the system
of first order ordinary differential equations given by Newtonian vector
field $\boldsymbol\Phi$ in tangent bundle $TM$:
$$
\hskip -2em
\boldsymbol\Phi=\sum^n_{i=1}v^i\cdot\frac{\partial}{\partial x^i}+
\sum^n_{i=1}\Phi^i\cdot\frac{\partial}{\partial v^i}.
\tag3.1
$$
Let $q=(p,\bold v)$ be a point of tangent bundle $TM$. Then canonical
projection $\pi\!:\nobreak TM\to M$ maps it to the point $p\in M$.
Canonical projection induces linear map $\pi_*\!: T_q(TM)
\to T_p(M)$. Applying this map to vector \thetag{3.1} at the point
$q\in TM$, we get velocity vector $\bold v$ at the point $p\in M$:
$$
\hskip -2em
\pi_*(\boldsymbol\Phi)=\bold v.
\tag3.2
$$
The equality \thetag{3.2} can be taken for the definition of Newtonian
vector field \thetag{3.1}. Linear map $\pi_*$ in \thetag{3.2} acts from
$T_q(TM)$ to $T_p(M)$. However, one can consider a map acting in opposite
direction.
\definition{Definition 3.1} Suppose that for each point $q$ of tangent
bundle $TM$ we have linear map $f\!: T_p(M)\to T_q(TM)$ from tangent
space $T_p(M)$ at the point $p=\pi(q)$ to tangent space $T_q(TM)$ at
the point $q$. This construction is called a {\bf lift} of vectors
from $M$ to tangent bundle $TM$.
\enddefinition
\definition{Definition 3.2} Lift of vectors $f$ from $M$ to $TM$ is
called {\bf smooth lift} if it maps each smooth vector field in $M$
to a smooth vector field in $TM$.
\enddefinition
\definition{Definition 3.3} Lift of vectors $f$ from $M$ to $TM$
is called {\bf vertical lift} if composition $\pi_*\compos f$ is
identically zero: $\pi_*\compos f=0$.
\enddefinition
\definition{Definition 3.4} Lift of vectors $f$ from $M$ to $TM$
is called {\bf horizontal lift} if composition $\pi_*\compos f$
is the field of identical operators on $M$, i\.~e\. $\pi_*\compos
f=\id$.
\enddefinition
   Each smooth manifold $M$ possesses canonical vertical lift of vectors
from $M$ to $TM$. It is defined as follows. Let $\bold X$ be some vector
field in $M$ and let $\bold X_p$ be its value at the point $p\in M$. Then
we can define one-parametric group of diffeomorphisms $\varphi_t$ in $TM$
that maps point $q=(p,\bold v)$ to the point $\varphi_t(q)=(p,\bold v+t
\cdot\bold X_p)$. \pagebreak Vector field $\bold Y$ in tangent bundle $TM$
associated with one-parametric group of diffeomorphisms $\varphi_t$ is
taken for the result of lifting vector field $\bold X$:
$$
\hskip -2em
\bold Y=w(\bold X).
\tag3.3
$$
Let's write \thetag{3.3} in local chart. Suppose that $\bold X$ is
given by its coordinates:
$$
\hskip -2em
\bold X=X^1\cdot\frac{\partial}{\partial x^1}+\ldots+X^n\cdot
\frac{\partial}{\partial x^n}.
\tag3.4
$$
Then result of applying canonical vertical lift $w$ to $\bold X$
is given by formula
$$
\hskip -2em
\bold Y=w(\bold X)=X^1\cdot\frac{\partial}{\partial v^1}+\ldots
+X^n\cdot\frac{\partial}{\partial v^n}.
\tag3.5
$$
Looking at \thetag{3.4} and \thetag{3.5}, we see that for each point
$q$ of tangent bundle $w$ is injective map that maps tangent space
$T_p(M)$ onto the vertical subspace
$$
\hskip -2em
V_q(TM)=\Ker\pi_*.
\tag3.6
$$
Let $h$ be some horizontal lift of vectors from $M$ to $TM$. For each
point $q$ of tangent bundle $h$ is also an injective. Let's denote by
$H_q(TM)$ its image
$$
\hskip -2em
H_q(TM)=\Img h.
\tag3.7
$$
Subspaces \thetag{3.6} and \thetag{3.7} are transversal and complementary
to each other:
$$
\hskip -2em
T_q(TM)=H_q(TM)\oplus V_q(TM).
\tag3.8
$$
\proclaim{Lemma 3.1} Defining horizontal lift of vectors from $M$ to $TM$
is equivalent to fixing horizontal subspace $H_q(TM)$ complementary to
vertical subspace $V_q(TM)$ in $T_q(TM)$ at each point $q$ of tangent
bundle $TM$.
\endproclaim
    Suppose that we have some horizontal lift $h$. Let's denote by
$\bold E_1,\,\ldots,\,\bold E_n$ base of coordinate vector fields
for some local chart of the manifold $M$:
$$
\hskip -2em
\bold E_1=\frac{\partial}{\partial x^1},\ \ldots,\ \bold E_n
=\frac{\partial}{\partial x^n}.
\tag3.9
$$
Due to \thetag{3.8} and due to the equality $\pi_*\compos h=\id$,
applying $h$ to $\bold E_i$, we get:
$$
\hskip -2em
h(\bold E_i)=\frac{\partial}{\partial x^i}-\sum^n_{k=1}
\Gamma^k_i\cdot\frac{\partial}{\partial v^k}.
\tag3.10
$$
Quantities $\Gamma^k_i=\Gamma^k_i(x^1,\ldots,x^n,v^1,\ldots,v^n)$
represent horizontal lift $h$ in local chart. If they are smooth
functions of their arguments, then $h$ is smooth lift. Under the
change of local chart in $M$ these quantities are transformed as
follows:
$$
\hskip -2em
\Gamma^k_i=\sum^n_{m=1}\sum^n_{a=1}S^k_m\,T^a_i\,\tilde\Gamma^m_{a}
+\sum^n_{m=1}\sum^n_{r=1}S^k_m\,\frac{\partial T^m_r}{\partial x^i}
\,v^r.
\tag3.11
$$
Here $S^k_m$, $T^a_i$, and $T^m_r$ are components of direct and inverse
transition matrices (Jacobi matrices) for the change of local coordinates:
$$
\xalignat 2
&\hskip -2em
S^i_j=\frac{\partial x^i}{\partial \tx^j},
&&T^i_j=\frac{\partial \tx^i}{\partial x^j}.
\tag3.12
\endxalignat
$$\par
   In arbitrary smooth manifold there is no canonical horizontal lift
of vectors. However, if Newtonian dynamical system \thetag{1.16} is
given, then there is horizontal lift $h$, canonically associated with
it (see \cite{22} and references therein). Let $\boldsymbol\Phi$ be
Newtonian vector field \thetag{3.1} and let $\bold X$ be some arbitrary
vector field in $TM$. Then we associate with $\bold X$ the following
vector field denoted by $\bold H(\bold X)$:
$$
\hskip -2em
\bold H(\bold X)=\frac{\bold X+[w\compos\pi_*(\bold X),\,
\boldsymbol\Phi]-w\compos\pi_*([\bold X,\,\boldsymbol\Phi])}{2}.
\tag3.13
$$
Here by square brackets we denote commutator of two vector fields.
\proclaim{Lemma 3.2} For any smooth function $\varphi$ and for arbitrary
smooth vector field $\bold X$ in $TM$ the equality $\bold H(\varphi\cdot
\bold X)=\varphi\cdot\bold H(\bold X)$ is fulfilled.
\endproclaim
\demo{Proof} By direct calculations we find that $[\varphi\cdot\bold X,
\,\boldsymbol\Phi]=\varphi\cdot[\bold X,\,\boldsymbol\Phi]-\boldsymbol
\Phi\varphi\cdot\bold X$. Here by $\boldsymbol\Phi\varphi$ we denote
derivative of function $\varphi$ along vector $\boldsymbol\Phi$:
$$
\boldsymbol\Phi\varphi=\sum^n_{i=1}v^i\,\frac{\partial\varphi}
{\partial x^i}+\sum^n_{i=1}\Phi^i\,\frac{\partial\varphi}
{\partial v^i}.
$$
Composition $w\compos\pi_*$ acts as linear operator at each point
$q\in TM$. Therefore
$$
\allowdisplaybreaks
\align
\hskip -2em
&w\compos\pi_*([\varphi\cdot\bold X,\,\boldsymbol\Phi])=\varphi
\cdot w\compos\pi_*([\bold X,\,\boldsymbol\Phi])-\boldsymbol\Phi\varphi
\cdot w\compos\pi_*(\bold X),
\tag3.14\\
\vspace{1ex}
\hskip -2em
&w\compos\pi_*(\varphi\cdot\bold X)=\varphi\cdot w\compos\pi_*(\bold X),
\\
\vspace{1ex}
\hskip -2em
&[w\compos\pi_*(\varphi\cdot\bold X),\,\boldsymbol\Phi]=\varphi\cdot
[w\compos\pi_*(\bold X),\,\boldsymbol\Phi]-\boldsymbol\Phi\varphi
\cdot w\compos\pi_*(\bold X).
\tag3.15
\endalign
$$
When subtracting \thetag{3.14} from \thetag{3.15}, terms containing
$\boldsymbol\Phi\varphi$ cancel each other. This proves required
equality $\bold H(\varphi\cdot\bold X)=\varphi\cdot\bold H(\bold X)$.
\qed\enddemo
    Formula \thetag{3.13} defines an operator $\bold H$ acting on
vector field $\bold X$ and yielding another vector field $\bold H
(\bold X)$. Commutators in \thetag{3.13} contain differentiation
with respect to $\bold X$. Therefore one might expect $\bold H$ to
be differential operator. Lemma~3.2 means that $\bold H$ is
non-differential linear operator acting pointwise at each point $q$
of tangent bundle $TM$.
\proclaim{Lemma 3.3} Linear operator $\bold H\!:T_q(TM)\to T_q(TM)$
defined by formula \thetag{3.13} satisfies the equalities $\bold H
\compos\bold H=\bold H$, $\pi_*\compos\bold H=\pi_*$, and $\bold H
\compos w=0$.
\endproclaim
\demo{Proof} In order to prove required equalities we shall use
direct calculations in local coordinates. Let's take an arbitrary
vector field $\bold X$ in $TM$ represented as
$$
\hskip -2em
\bold X=\sum^n_{i=1}X^i\cdot\frac{\partial}{\partial x^i}
+\sum^n_{i=1}Y^i\cdot\frac{\partial}{\partial v^i}.
\tag3.16
$$
Newtonian vector field $\boldsymbol\Phi$ is represented by formula
\thetag{3.1}. Hence for $[\bold X,\,\boldsymbol\Phi]$ we have
$$
\gather
[\bold X,\,\boldsymbol\Phi]=\sum^n_{k=1}Y^k\cdot\frac{\partial}
{\partial x^k}+\sum^n_{k=1}\sum^n_{i=1}\left(\!X^i\,\frac{\partial
\Phi^k}{\partial x^i}+Y^i\,\frac{\partial\Phi^k}{\partial v^i}
\right)\cdot\frac{\partial}{\partial v^k}\,-\\
\vspace{1ex}
-\sum^n_{i=1}\sum^n_{k=1}\left(\!v^k\,\frac{\partial X^i}{\partial x^k}
+\Phi^k\,\frac{\partial X^i}{\partial v^k}\right)\cdot\frac{\partial}
{\partial x^i}
-\sum^n_{i=1}\sum^n_{k=1}\left(\!v^k\,\frac{\partial Y^i}{\partial x^k}
+\Phi^k\,\frac{\partial Y^i}{\partial v^k}\right)\cdot\frac{\partial}
{\partial v^i}.
\endgather
$$
Applying operator $w\compos\pi_*$ to \thetag{3.16} and to the above
expression, we obtain
$$
\gather
\hskip -6em
w\compos\pi_*(\bold X)=\sum^n_{i=1}X^i\cdot\frac{\partial}
{\partial v^i},
\tag3.17
\\
\vspace{1ex}
\hskip -6em
w\compos\pi_*([\bold X,\,\boldsymbol\Phi])=\sum^n_{k=1}Y^k\cdot
\frac{\partial}{\partial v^k}
-\sum^n_{i=1}\sum^n_{k=1}\left(\!v^k\,\frac{\partial X^i}{\partial x^k}
+\Phi^k\,\frac{\partial X^i}{\partial v^k}\right)\cdot\frac{\partial}
{\partial v^i}.\hskip -2em
\tag3.18
\endgather
$$
Now, relying upon formula \thetag{3.17}, we calculate the commutator
$[w\compos\pi_*(\bold X),\,\boldsymbol\Phi]$:
$$
\hskip -2em
\gathered
[w\compos\pi_*(\bold X),\,\boldsymbol\Phi]=\sum^n_{k=1}X^k\cdot
\frac{\partial}{\partial x^k}+\sum^n_{k=1}\sum^n_{i=1}
X^i\,\frac{\partial\Phi^k}{\partial v^i}\cdot\frac{\partial}
{\partial v^k}\,-\\
\vspace{1ex}
-\sum^n_{i=1}\sum^n_{k=1}\left(\!v^k\,\frac{\partial X^i}{\partial x^k}
+\Phi^k\,\frac{\partial X^i}{\partial v^k}\right)\cdot\frac{\partial}
{\partial v^i}.
\endgathered
\tag3.19
$$
Then subtract \thetag{3.18} from \thetag{3.19}. As a result we
obtain the following equality:
$$
\gather
[w\compos\pi_*(\bold X),\,\boldsymbol\Phi]-w\compos\pi_*([\bold X,\,
\boldsymbol\Phi])=\sum^n_{k=1}X^k\cdot\frac{\partial}{\partial x^k}\,+\\
\vspace{1ex}
+\sum^n_{k=1}\sum^n_{i=1}X^i\,\frac{\partial\Phi^k}{\partial v^i}
\cdot\frac{\partial}{\partial v^k}-\sum^n_{k=1}Y^k\cdot
\frac{\partial}{\partial v^k}.
\endgather
$$
And ultimately, for the result of applying operator $\bold H$ to 
vector $\bold X$ given by formula \thetag{3.16} we obtain the
following expression:
$$
\hskip -2em
\bold H(\bold X)=\sum^n_{k=1}X^k\cdot\frac{\partial}{\partial x^k}
+\frac{1}{2}\sum^n_{s=1}\sum^n_{i=1}X^i\,\frac{\partial\Phi^s}
{\partial v^i}\cdot\frac{\partial}{\partial v^s}.
\tag3.20
$$
Looking at the expression in right hand side of \thetag{3.20},
we see that the required equalities $\bold H\compos\bold H=\bold H$,
$\pi_*\compos\bold H=\pi_*$, and $\bold H\compos w=0$ for operator
$\bold H$ appear to be obvious. Lemma~3.3 is proved.
\qed\enddemo
    First equality $\bold H\compos\bold H=\bold H$ means that $\bold H$
is an operator of projection. It projects $T_q(TM)$ onto some subspace
$H_q(TM)$, where $H_q(TM)=\Img\bold H$. Each projection operator breaks
the space into direct sum of its kernel and its image: 
$$
\hskip -2em
T_q(TM)=\Ker\bold H\oplus\Img\bold H.
\tag3.21
$$
Second equality $\pi_*\compos\bold H=\pi_*$ yields $\Ker\bold H\subseteq
\Ker\pi_*$, while third equality $\bold H\compos w=0$ means that $\Img w
\subseteq\Ker\bold H$. Taking into account \thetag{3.6}, we obtain that
kernel of operator $\bold H$ coincides with vertical subspace $V_q(TM)$:
$$
\hskip -2em
\Ker\bold H=\Ker\pi_*=\Img w=V_q(TM).
\tag3.22
$$
Due to \thetag{3.22} the equality \thetag{3.21} can be rewritten as
\thetag{3.8}. This means that operator $\bold H$ defines horizontal
subspace $H_q(TM)=\Img\bold H$ in $T_q(TM)$. Applying lemma~3.1, we
derive the following theorem.
\proclaim{Theorem 3.1} Each Newtonian dynamical system \thetag{1.16}
in smooth manifold $M$ generates horizontal lift of vectors from $M$
to $TM$ canonically associated with it.
\endproclaim
    This result is not new (see \cite{22} and references therein). We
reproduced it here for the sake of completeness of our consideration.
\par
    Let $h$ be horizontal lift canonically associated with Newtonian
dynamical system \thetag{1.16}. Then $\bold H=h\compos\pi_*$ and
$h(\bold E_i)=h\compos\pi_*(\partial/\partial x^i)=\bold H(\partial/
\partial x^i)$. Applying formula \thetag{3.20} to vector field $\bold
X=\bold E_i$, we now obtain the following equality:
$$
\hskip -2em
h(\bold E_i)=\bold H(\partial/\partial x^i)=\frac{\partial}
{\partial x^i}+\sum^n_{k=1}\frac{1}{2}\,\frac{\partial\Phi^k}
{\partial v^i}\cdot\frac{\partial}{\partial v^k}.
\tag3.23
$$
Comparing \thetag{3.23} with \thetag{3.10}, for components of
$h$ we derive explicit formula:
$$
\hskip -2em
\Gamma^k_i=-\frac{1}{2}\,\frac{\partial\Phi^k}{\partial v^i}.
\tag3.24
$$
\head
4. Extended tensor fields and covariant differentiations.
\endhead
    Concept of extended tensor fields is closely related to
Newtonian dynamical system. Indeed, if we look at right hand
side of \thetag{1.4}, we see that $\bold F=\bold F(p,\bold v)$
is a vector-valued function with values in tangent spaces $T_p(M)$.
However, it doesn't fit standard definition of vector field in
$M$ since it depends not only on a point $p\in M$, but also
upon velocity vector $\bold v$ at this point. Both $p$ and $\bold v$
form a point $q=(p,\bold v)$ of tangent bundle $TM$. Function
$\bold F=\bold F(p,\bold v)$ is an example of extended vector field.
Vector-function $\bold V=\bold V(p,\bold p)$ in \thetag{2.3} with
components \thetag{2.4} is another example. Its arguments form a
point $q=(p,\bold p)$ if cotangent bundle $T^*\!M$.\par
    In order to formulate definition of extended tensor fields
in general let's consider the following tensor product of tangent
spaces and their dual spaced:
$$
T^r_s(p,M)=\overbrace{T_p(M)\otimes\ldots\otimes T_p(M)}^{\text{$r$
times}}\otimes\underbrace{T^*_p(M)\otimes\ldots\otimes T^*_p(M)}_{\text{$s$
times}}
$$
\definition{Definition 4.1}Extended tensor field $\bold X$ of type
$(r,s)$ in $\bold v$-representation is a tensor-valued function with
argument $q=(p,\bold v)$ in tangent bundle $TM$ and with value
$\bold X(q)$ in tensor space $T^r_s(p,M)$, where $p=\pi(q)$.
\enddefinition
\definition{Definition 4.2}Extended tensor field $\bold X$ of type
$(r,s)$ in $\bold p$-representation is a tensor-valued function with
argument $q=(p,\bold p)$ in cotangent bundle $T^*\!M$ and with value
$\bold X(q)$ in tensor space $T^r_s(p,M)$, where $p=\pi(q)$.
\enddefinition
    At first let's consider extended tensor fields in
$\bold v$-representation. Denote by $T^r_s(M)$ the set of smooth
extended tensor fields of type $(r,s)$ and take the sum
$$
\hskip -2em
\bold T(M)=\bigoplus^\infty_{r=0}\bigoplus^\infty_{s=0}T^r_s(M)
\tag4.1
$$
In $\bold T(M)$ we have all standard tensorial operations like summation,
multiplication by scalars, tensor product, and contraction. Direct sum
\thetag{4.1} possesses structure of graded algebra over the ring of smooth
scalar functions in $TM$, we denote this ring by $\goth F(TM)$. Algebra
$\bold T(M)$ is called {\bf extended algebra of tensor fields}.
\definition{Definition 4.3} A map $D\!:\bold T(M)\to\bold T(M)$ is
called {\bf differentiation} of extended algebra of tensor fields,
if the following conditions are fulfilled:
\roster
\rosteritemwd=10pt
\item concordance with grading: $D(T^r_s(M))\subset T^r_s(M)$;
\item $\Bbb R$-linearity: $D(\bold X+\bold Y)=D(\bold X)+D(\bold Y)$
      and $D(\lambda\bold X)=\lambda D(\bold X)$ for $\lambda\in\Bbb R$;
\item commutation with contractions: $D(C(\bold X))=
      C(D(\bold X))$;
\item Leibniz rule: $D(\bold X\otimes\bold Y)=D(\bold X)
      \otimes\bold Y+\bold X\otimes D(\bold Y)$.
\endroster
\enddefinition
     Let's denote by $\goth D(M)$ the total set of differentiations
of extended algebra of tensor fields $\bold T(M)$. It is easy to see
that it possesses the structure of $\goth F(TM)$-module. The set of
extended vector fields $T^1_0(M)$ is equipped with the same structure
of module over the ring $\goth F(TM)$. This coincidence motivates
the following definition.
\definition{Definition 4.4} {\bf Covariant differentiation} $\nabla$ in the
algebra of extended tensor fields $\bold T(M)$ is a homomorphism of $\goth
F(TM)$-modules $\nabla\!:T^1_0(M)\to\goth D(M)$. Image of vector field
$\bold Y$ under such homomorphism is called covariant differentiation
along vector field $\bold Y$. It is denoted by $\nabla_{\bold Y}$.
\enddefinition
   In Chapter~\uppercase\expandafter{\romannumeral 3} of thesis \cite{7} it
was shown that each covariant differentiation $\nabla$ in extended algebra
of tensor fields is associated with some lift of vectors from $M$ to $TM$.
It is called {\bf horizontal covariant differentiation} (or {\bf vertical
covariant differentiation}) if corresponding lift of vectors is horizontal
(or vertical). In each smooth manifold there is canonical vertical
covariant differentiation associated with canonical vertical lift of
vectors $w$. We denote it by $\tilde\nabla$. In local chart this covariant
differentiation is given by the following formula:
$$
\hskip -2em
\tilde\nabla_mX^{i_1\ldots\,i_r}_{j_1\ldots\,j_s}=\frac{\partial
X^{i_1\ldots\,i_r}_{j_1\ldots\,j_s}}{\partial v^m}.
\tag4.2
$$
Due to \thetag{4.2} this covariant differentiation is also called
{\bf velocity gradient}.\par
    To define horizontal covariant differentiation, apart from horizontal
lift of vectors, one need some extended affine connection $\Gamma$. Its
components are given by formula
$$
\hskip -2em
\nabla_{\bold E_{\ssize i}}\bold E_j=\sum^n_{k=1}
\Gamma^k_{ij}\,\bold E_k.
\tag4.3
$$
Here $\bold E_1,\,\ldots,\,\bold E_n$ are coordinate vector fields
\thetag{3.9} in $M$. Unlike traditional affine connection, components
of extended affine connection $\Gamma^k_{ij}$ depend on double set of
arguments, i\.\,e\. on coordinates of point $p\in M$ and on components
of vector $\bold v$:
$$
\hskip -2em
\Gamma^k_{ij}=\Gamma^k_{ij}(x^1,\ldots,x^n,v^1,\ldots,v^n).
\tag4.4
$$
Under the change of local chart quantities \thetag{4.4} are transformed
as follows:
$$
\hskip -2em
\Gamma^k_{ij}=\sum^n_{m=1}\sum^n_{a=1}\sum^n_{c=1} S^k_m\,T^a_i
\,T^c_j\,\tilde\Gamma^m_{ac}+\sum^n_{m=1} S^k_m\,\frac{\partial
T^m_i}{\partial x^j}
\tag4.5
$$
(matrices $S$ and $T$ are defined in \thetag{3.12}).
If components of horizontal lift $h$ in \thetag{3.10} and components
of extended affine connection $\Gamma$ in \thetag{4.3} are given, then
horizontal covariant differentiation $\nabla$ in local coordinates is
given by the following formula:
$$
\hskip -2em
\aligned
&\nabla_mX^{i_1\ldots\,i_r}_{j_1\ldots\,j_s}=\frac{\partial
X^{i_1\ldots\,i_r}_{j_1\ldots\,j_s}}{\partial x^m}
-\sum^n_{a=1}\sum^n_{b=1}\Gamma^b_m\,\frac{\partial
X^{i_1\ldots\,i_r}_{j_1\ldots\,j_s}}{\partial v^b}\,+\\
&+\sum^r_{k=1}\sum^n_{a_k=1}\!\Gamma^{i_k}_{m\,a_k}\,X^{i_1\ldots\,
a_k\ldots\,i_r}_{j_1\ldots\,\ldots\,\ldots\,j_s}
-\sum^s_{k=1}\sum^n_{b_k=1}\!\Gamma^{b_k}_{m\,j_k}
X^{i_1\ldots\,\ldots\,\ldots\,i_r}_{j_1\ldots\,b_k\ldots\,j_s}.
\endaligned
\tag4.6
$$\par
Note that $\Gamma^k_i$ and $\Gamma^k_{ij}$ are two independent sets
of parameters in \thetag{4.6}. The only condition is that they should
obey transformation rules \thetag{3.11} and \thetag{4.5} respectively.
However, one can prove the following two lemmas. 
\proclaim{Lemma 4.1} If horizontal lift of vectors from $M$ to $TM$ is
given and if $\Gamma^k_i$ are its components, then there is symmetric
extended affine connection with components
$$
\hskip -2em
\Gamma^k_{ij}=\frac{1}{2}\,\frac{\partial\Gamma^k_i}{\partial v^j}
+\frac{1}{2}\,\frac{\partial\Gamma^k_j}{\partial v^i}.
\tag4.7
$$
\endproclaim
\proclaim{Lemma 4.2} If extended affine connection with components
$\Gamma^k_{ij}$ is given, then there is horizontal lift of vectors from
$M$ to $TM$ with components
$$
\hskip -2em
\Gamma^k_i=\sum^n_{j=1}\Gamma^k_{ij}\,v^j.
\tag4.8
$$
\endproclaim
These two lemmas are proved by direct calculations on the base of
formulas \thetag{3.11} and \thetag{4.5}. Applying \thetag{4.7} to
\thetag{3.24}, we obtain
$$
\hskip -2em
\Gamma^k_{ij}=-\frac{1}{2}\frac{\partial^2\Phi^k}{\partial v^i\,
\partial v^j}.
\tag4.9
$$
As a result we have proved the following theorem.
\proclaim{Theorem 4.1} Each Newtonian dynamical system \thetag{1.16} in
smooth manifold $M$ generates extended affine connection with components
\thetag{4.9} canonically associated with this dynamical system.
\endproclaim
    Usually equalities \thetag{4.7} and \thetag{4.8} cannot be fulfilled
simultaneously. In all previous papers we defined $\Gamma^k_i$ through
$\Gamma^k_{ij}$ by means of formula \thetag{4.8}. This is equivalent to
the equality $\nabla_iv^k=0$. In present paper we keep this tradition.
Therefore formula \thetag{4.6} for horizontal covariant derivative looks
like
$$
\hskip -2em
\aligned
&\nabla_mX^{i_1\ldots\,i_r}_{j_1\ldots\,j_s}=\frac{\partial
X^{i_1\ldots\,i_r}_{j_1\ldots\,j_s}}{\partial x^m}
-\sum^n_{a=1}\sum^n_{b=1}\sum^n_{c=1}v^c\,\Gamma^b_{cm}\,\frac{\partial
X^{i_1\ldots\,i_r}_{j_1\ldots\,j_s}}{\partial v^b}\,+\\
&+\sum^r_{k=1}\sum^n_{a_k=1}\!\Gamma^{i_k}_{m\,a_k}\,X^{i_1\ldots\,
a_k\ldots\,i_r}_{j_1\ldots\,\ldots\,\ldots\,j_s}
-\sum^s_{k=1}\sum^n_{b_k=1}\!\Gamma^{b_k}_{m\,j_k}
X^{i_1\ldots\,\ldots\,\ldots\,i_r}_{j_1\ldots\,b_k\ldots\,j_s}.
\endaligned
\tag4.10
$$\par
    Theory of differentiation for extended tensor fields in
$\bold p$-representation is a little bit different. Here we also can
define {\bf extended algebra of tensor fields } by means of direct
sum \thetag{4.1}. It is algebra over the ring of smooth functions
in cotangent bundle $\goth F(T^*\!M)$. Definition~4.3 remains unchanged.
Definition~4.4 is replaced by the following two definitions.
\definition{Definition 4.5} {\bf Covariant differentiation} $\nabla$ in the
algebra of extended vector fields $\bold T(M)$ is a homomorphism of $\goth
F(T^*\!M)$-modules $\nabla\!:T^1_0(M)\to\goth D(M)$. Image of vector field
$\bold Y$ under such homomorphism is called covariant differentiation
along vector field $\bold Y$. It is denoted by $\nabla_{\bold Y}$.
\enddefinition
\definition{Definition 4.6} {\bf Contravariant differentiation} $\nabla$
in the algebra of extended vector fields $\bold T(M)$ is a homomorphism
of $\goth F(T^*\!M)$-modules $\nabla\!:T^0_1(M)\to\goth D(M)$. Image of
covector field $\bold q$ under such homomorphism is called contravariant
differentiation along covector field $\bold q$. It is denoted by
$\nabla_{\bold q}$.
\enddefinition
   Instead of velocity gradient \thetag{4.2} here in
$\bold p$-representation we have {\bf canonical vertical contravariant
differentiation} $\tilde\nabla$. It is given by formula
$$
\hskip -2em
\tilde\nabla^mX^{i_1\ldots\,i_r}_{j_1\ldots\,j_s}=\frac{\partial
X^{i_1\ldots\,i_r}_{j_1\ldots\,j_s}}{\partial p_m}.
\tag4.11
$$
Due to \thetag{4.11} contravariant differentiation $\tilde\nabla$ is
called {\bf momentum gradient}. Instead of canonical vertical lift of
vectors in $\bold p$-representation for each point $q=(p,\bold p)$ of
cotangent bundle $T^*\!M$ we have injective map $w\!: T^*_p(p)\to
T_q(TM)$ associated with differentiation $\tilde\nabla$. Suppose that
covector $\bold q$ is given by its coordinates:
$$
\hskip -2em
\bold q=q_1\cdot dx^1+\ldots+q_n\cdot dx^n.
\tag4.12
$$
Applying $w$ to $\bold q$, we obtain a vector $\bold Y=w(\bold q)$
given by the following expression:
$$
\hskip -2em
\bold Y=w(\bold q)=q_1\cdot\frac{\partial}{\partial p_1}+\ldots
+q_n\cdot\frac{\partial}{\partial p_n}.
\tag4.13
$$
Formulas \thetag{4.12} and \thetag{4.13} are similar to formulas
\thetag{3.4} and \thetag{3.5}. However, there is invariant way to
determine vector $\bold Y$. It is generated by one-parametric group
of diffeomorphisms $\varphi_t$ in $T^*\!M$ which is defined as
follows. If $q=(p,\bold p)$ is a point of cotangent bundle $T^*\!M$,
then $\varphi_t(q)=(p,\bold p+t\cdot\bold q)$.\par
   Injective map $w\!: T^*_p(p)\to T_q(TM)$ defined just above satisfy
the equality $\pi_*\compos w=0$. Its image coincides with vertical
subspace in $T_q(T^*\!M)$:
$$
\Img w=V_q(T^*\!M)=\Ker\pi_*.
$$
Each horizontal covariant derivative in $\bold p$-representation
is given by formula
$$
\aligned
&\nabla_mX^{i_1\ldots\,i_r}_{j_1\ldots\,j_s}=\frac{\partial
X^{i_1\ldots\,i_r}_{j_1\ldots\,j_s}}{\partial x^m}
+\sum^n_{a=1}\sum^n_{b=1}\Gamma_{mb}\,\frac{\partial
X^{i_1\ldots\,i_r}_{j_1\ldots\,j_s}}{\partial p_b}\,+\\
&+\sum^r_{k=1}\sum^n_{a_k=1}\!\Gamma^{i_k}_{m\,a_k}\,X^{i_1\ldots\,
a_k\ldots\,i_r}_{j_1\ldots\,\ldots\,\ldots\,j_s}
-\sum^s_{k=1}\sum^n_{b_k=1}\!\Gamma^{b_k}_{m\,j_k}
X^{i_1\ldots\,\ldots\,\ldots\,i_r}_{j_1\ldots\,b_k\ldots\,j_s}.
\endaligned
$$
Here $\Gamma_{ij}$ and $\Gamma^k_{ij}$ are components of some
horizontal lift of vectors from $M$ to $T^*\!M$ and some extended
affine connection respectively. Quantities $\Gamma^k_{ij}$ are
defined by the equality \thetag{4.3}. They satisfy the same
transformation rule \thetag{4.5} as quantities \thetag{4.4}. However,
they differ from \thetag{4.4} since they have another set of arguments:
$$
\hskip -2em
\Gamma^k_{ij}=\Gamma^k_{ij}(x^1,\ldots,x^n,p_1,\ldots,p_n).
\tag4.14
$$
Horizontal lift of vectors $h$ from $M$ to $T^*\!M$ satisfies the
equality $\pi_*\compos h=\id$. It's components are determined by
the following formula:
$$
h(\bold E_i)=\frac{\partial}{\partial x^i}+\sum^n_{k=1}
\Gamma_{ik}\cdot\frac{\partial}{\partial p_k}.
$$
Under change of local chart quantities $\Gamma_{ij}$ are transformed
as follows:
$$
\hskip -2em
\Gamma_{ij}=\sum^n_{a=1}\sum^n_{c=1}T^a_i
\,T^c_j\,\tilde\Gamma_{ac}+\sum^n_{k=1}\sum^n_{m=1}
p_k\,S^k_m\,\frac{\partial
T^m_i}{\partial x^j}.
\tag4.15
$$
Comparing formula \thetag{4.15} with \thetag{4.5}, one can prove two
lemmas which are similar to lemma~4.1 and lemma~4.2.
\proclaim{Lemma 4.3} If horizontal lift of vectors from $M$ to $T^*\!M$
is given and if $\Gamma_{ij}$ are its components, then there is symmetric
extended affine connection with components
$$
\hskip -2em
\Gamma^k_{ij}=\frac{\partial\Gamma_{ij}}{\partial p_k}.
\tag4.16
$$
\endproclaim
\proclaim{Lemma 4.4} If extended affine connection with components
$\Gamma^k_{ij}$ is given, then there is horizontal lift of vectors from
$M$ to $T^*\!M$ with components
$$
\hskip -2em
\Gamma_{ij}=\sum^n_{k=1}\Gamma^k_{ij}\,p_k.
\tag4.17
$$
\endproclaim
    Usually equalities \thetag{4.16} and \thetag{4.17} cannot be fulfilled
simultaneously. In previous papers we defined $\Gamma_{ij}$ by means of
formula \thetag{4.17}. This is equivalent to
$$
\nabla_ip_j=0.
$$
In present paper we keep this tradition. Connection components
$\Gamma^k_{ij}$ will be imported to $\bold p$-representation by
means of generalized \pagebreak Legendre transformation \thetag{2.2}
(see below). Then formula for horizontal covariant differentiation
looks like
$$
\hskip -2em
\aligned
&\nabla_mX^{i_1\ldots\,i_r}_{j_1\ldots\,j_s}=\frac{\partial
X^{i_1\ldots\,i_r}_{j_1\ldots\,j_s}}{\partial x^m}
-\sum^n_{a=1}\sum^n_{b=1}\sum^n_{c=1}p_c\,\Gamma^c_{mb}\,\frac{\partial
X^{i_1\ldots\,i_r}_{j_1\ldots\,j_s}}{\partial p_b}\,+\\
&+\sum^r_{k=1}\sum^n_{a_k=1}\!\Gamma^{i_k}_{m\,a_k}\,X^{i_1\ldots\,
a_k\ldots\,i_r}_{j_1\ldots\,\ldots\,\ldots\,j_s}
-\sum^s_{k=1}\sum^n_{b_k=1}\!\Gamma^{b_k}_{m\,j_k}
X^{i_1\ldots\,\ldots\,\ldots\,i_r}_{j_1\ldots\,b_k\ldots\,j_s}.
\endaligned
\tag4.18
$$
Formulas \thetag{4.10} and \thetag{4.18} mean that we do not use
horizontal lift with components \thetag{3.24} directly, but only
for deriving extended connection with components \thetag{4.9}.
\head
5. Generalized Legendre transformation.
\endhead
   In section~2 we learned that some extended affine connection and some
horizontal lift of vectors from $M$ to $TM$ are canonically associated
with Newtonian dynamical system \thetag{1.16}. As for dynamical system
\thetag{2.3} in $T^*\!M$, I don't know similar direct constructions for
$\Gamma^k_{ij}$ and $\Gamma_{ij}$ in this case. Therefore now we shall
import these quantities from $\bold v$-representation to
$\bold p$-representation by means of generalized Legendre transformation
$\lambda$. For $\Gamma^k_{ij}$ in \thetag{4.14} we write:
$$
\hskip -2em
\Gamma^k_{ij}=\Gamma^k_{ij}\compos\lambda^{-1}.
\tag5.1
$$
This equality \thetag{5.1} means that we simply change arguments
in \thetag{4.4} by substituting functions \thetag{2.4} for $v^1,\,
\ldots,\,v^n$.\par
    Quantities $\Gamma^k_{ij}$ in $\bold v$-representation are given by
formula \thetag{4.9}. Using formula \thetag{5.1}, one can find their
counterparts $\Gamma^k_{ij}$ in $\bold p$-representation. However, it
would be natural to express them through functions $V^1,\,\ldots,\,V^n$
and $\Theta_1,\,\ldots,\,\Theta_n$ that determine dynamical system
\thetag{2.3}. Differential equations \thetag{2.3} describe the same
dynamics as differential equations \thetag{1.16}, but in other variables
$x^1,\,\ldots,\,x^n,\,p_1,\,\ldots,\,p_n$. Therefore we can calculate
time derivative of velocity vector $\bold v$ as follows:
$$
\hskip -2em
\dot v^k=\frac{dV^k}{dt}=\sum^n_{i=1}\frac{\partial V^k}{\partial x^i}\,
V^i+\sum^n_{i=1}\frac{\partial V^k}{\partial p_i}\,\Theta_i.
\tag5.2
$$
Comparing \thetag{5.2} with second part of equations \thetag{1.16}, we
obtain the equality
$$
\hskip -2em
\Phi^k\compos\lambda^{-1}=\sum^n_{i=1}\frac{\partial V^k}{\partial x^i}\,
V^i+\sum^n_{i=1}\frac{\partial V^k}{\partial p_i}\,\Theta_i.
\tag5.3
$$
Note that partial derivatives $\partial V^k/\partial p_r$ in \thetag{5.3}
form Jacobi matrix for transformation given by functions \thetag{2.4},
when $x^1,\,\ldots,\,x^n$ are fixed. Let's denote
$$
\hskip -2em
g^{ir}=\frac{\partial V^i}{\partial p_r}=\tilde\nabla^rV^i.
\tag5.4
$$
Formula \thetag{5.4} determines components of an extended tensor field
$\bold g$. Tensor $\bold g$ plays the role similar to that of dual metric
tensor in Riemannian geometry. \pagebreak It is non-degenerate: $\det
\bold g\neq 0$, but, in general, it is not symmetric: $g^{ij}\neq g^{ji}$.
By $g_{ij}$ we denote components of inverse matrix for \thetag{5.4},
i\.\,e\.
$$
\xalignat 2
\hskip -2em
&\sum^n_{j=1}g^{ij}\,g_{jk}=\delta^i_k,
&&\sum^n_{j=1}g_{ij}\,g^{jk}=\delta^k_i.
\tag5.5
\endxalignat 
$$
They form another extended tensor field, which traditionally is denoted
by the same symbol $\bold g$. Though being non-symmetric extended tensor
field, it is direct analog of metric tensor in Riemannian geometry.\par
     Let $f=f(x^1,\ldots,x^n,v^1,\ldots,v^n)$ be some function in
$\bold v$-representation for some local chart of $M$. This might be
component of extended tensor field, component of extended connection,
or coordinate representation of some other geometric object. Then
$f\compos\lambda^{-1}$ is its $\bold p$-representation. Using explicit
form \thetag{2.4} of the map $\lambda^{-1}$, we can derive the following
transformation rules for partial derivatives:
$$
\align
&\hskip -2em
\frac{\partial(f\compos\lambda^{-1})}{\partial p_k}=\sum^n_{i=1}
\frac{\partial V^i}{\partial p_k}\cdot\left(\frac{\partial f}
{\partial v^i}\compos\lambda^{-1}\!\right)\!,
\tag5.6\\
\vspace{1ex}
&\hskip -2em
\frac{\partial(f\compos\lambda^{-1})}{\partial x^k}=\sum^n_{i=1}
\frac{\partial V^i}{\partial x^k}\cdot\left(\frac{\partial f}
{\partial v^i}\compos\lambda^{-1}\!\right)+\frac{\partial f}{\partial v^i}
\compos\lambda^{-1}.
\endalign
$$
Taking into account \thetag{5.4} and \thetag{5.5}, we can rewrite
\thetag{5.6} as
$$
\hskip -2em
\frac{\partial f}{\partial v^i}\compos\lambda^{-1}=
\sum^n_{k=1}g_{ki}\cdot\frac{\partial(f\compos\lambda^{-1})}
{\partial p_k}.
\tag5.7
$$
Now let's apply \thetag{5.7} to the function $f=\Phi^k(x^1,\ldots,x^n,
v^1,\ldots,v^n)$ in \thetag{5.3}:
$$
\frac{\partial\Phi^k}{\partial v^i}\compos\lambda^{-1}=
\sum^n_{r=1}g_{ri}\cdot\frac{\partial}{\partial p_r}
\left(\,\shave{\sum^n_{m=1}}\frac{\partial V^k}{\partial x^m}\,
V^m+\shave{\sum^n_{m=1}}\frac{\partial V^k}{\partial p_m}\,\Theta_m
\right)\!.
$$
Applying \thetag{5.7} once more and taking into account \thetag{4.9}
and \thetag{5.1}, we obtain
$$
\Gamma^k_{ij}=
\sum^n_{r=1}\sum^n_{s=1}g_{sj}\cdot\frac{\partial}{\partial p_s}\!
\left(g_{ri}\cdot\frac{\partial}{\partial p_r}\!
\left(\,\shave{\sum^n_{m=1}}\frac{\partial V^k}{\partial x^m}\,
V^m+\shave{\sum^n_{m=1}}\frac{\partial V^k}{\partial p_m}\,\Theta_m
\right)\!\right)\!.
$$
For the sake of convenience this formula should be slightly
transformed:
$$
\hskip -2em
\gathered
\Gamma^k_{ij}=\sum^n_{r=1}\sum^n_{s=1}g_{ri}\,g_{sj}\cdot
\frac{\partial^2}{\partial p_r\,\partial p_s}\left(\,\shave{\sum^n_{m=1}}
\frac{\partial V^k}{\partial x^m}\,V^m+\shave{\sum^n_{m=1}}
\frac{\partial V^k}{\partial p_m}\,\Theta_m\right)\,-\\
-\sum^n_{r=1}\sum^n_{s=1}\sum^n_{\alpha=1}g_{ri}\,g_{sj}\,
\frac{\partial^2V^\alpha}{\partial p_r\,\partial p_s}\cdot
\frac{\partial}{\partial p_\alpha}\!\left(\,\shave{\sum^n_{m=1}}
\frac{\partial V^k}{\partial x^m}\,V^m+\shave{\sum^n_{m=1}}
\frac{\partial V^k}{\partial p_m}\,\Theta_m\right)\!.
\endgathered
\tag5.8
$$
Looking at \thetag{5.8}, it is obvious that $\Gamma^k_{ij}$ keep
symmetry $\Gamma^k_{ij}=\Gamma^k_{ji}$ in $\bold p$-representation
as well. It is not surprising due to \thetag{5.1}, \pagebreak but
now it is explicitly evident.
\proclaim{Theorem 5.1} Each Newtonian dynamical system \thetag{2.3}
in cotangent bundle $T^*\!M$ generates extended affine connection with
components \thetag{5.8} canonically associated with this dynamical
system.
\endproclaim
    Theorem~5.1 reformulates previous theorem~4.1 with respect to
the same Newtonian dynamics, but transferred from tangent bundle
to cotangent bundle by means of generalized Legendre transformation
\thetag{2.2}.
\head
6. Regularity condition.
\endhead
   Now we are ready to study Newtonian dynamics given by the
equations \thetag{2.3} without referring to its $\bold v$-representation
\thetag{1.16}. If we consider trajectory of dynamical system as a curve
$p=p(t)$ in $M$, then velocity vector on this trajectory is given by
the value of extended vector field $\bold V=\bold V(p,\bold p)$. Let's
define another extended vector field $\bold W=\bold W(p,\bold p)$ with
components
$$
\hskip -2em
W^s=\sum^n_{r=1}\tilde\nabla^sV^r\,p_r
\tag6.1
$$
and consider scalar product of this vector field $\bold W$ and momentum
covector $\bold p$:
$$
\hskip -2em
\Omega=\left<\bold p\,|\,\bold W\right>=\sum^n_{s=1}p_s\,W^s(x^1,\ldots,
x^n,v^1,\ldots,v^n).
\tag6.2
$$
Extended scalar field $\Omega=\Omega(p,\bold p)$ in \thetag{6.2} is
somewhat like {\it kinetic energy} in mechanics (compare \thetag{6.2}
with \thetag{1.11} and \thetag{1.13}). In mechanics both momentum and
kinetic energy represent the ``amount of motion'' stored in moving
object. This motivates the following definition.
\definition{Definition 6.1} Generalized Legendre map \thetag{2.2}
given by components of extended vector field $\bold V$ in \thetag{2.4}
is called {\bf regular} if
\roster
\rosteritemwd=10pt
\item it is diffeomorphic;
\item $V(p,\bold p)=0$ is equivalent to $\bold p=0$;
\item $\Omega(p,\bold p)\neq 0$ for $\bold p\neq 0$.
\endroster
\enddefinition
    Third part of regularity condition means that vector $\bold W
=\bold W(p,\bold p)$ is transversal to null-space of momentum covector
$\bold p$. Therefore one can consider operator-valued extended tensor
field $\bold P$ with the following components:
$$
\hskip -2em
P^i_j=\delta^i_j-\frac{W^i\,p_j}{\Omega}.
\tag6.3
$$
Operator $\bold P=\bold P(p,\bold p)$ with components \thetag{6.3} is a
projector onto null-space of momentum covector $\bold p$ along vector
$\bold W=\bold W(p,\bold p)$. It is defined everywhere in $T^*\!M$
except for those points $q=(p,\bold p)$, where $\bold p=0$.
\head
7. Force covector.
\endhead
    Let's consider Newtonian dynamical system \thetag{2.3} in smooth
manifold $M$. Suppose that $M$ is equipped with symmetric extended
affine connection $\Gamma$. This might be connection \thetag{5.8}
canonically associated with dynamical system \thetag{2.3} or any other
symmetric extended affine connection which is not related to \thetag{2.3}
at all. In both cases one can introduce extended covector field
$\bold Q$ with components
$$
\hskip -2em
Q_i=\Theta_i-\sum^n_{j=1}\sum^n_{k=1}\Gamma^k_{ij}\,p_k\,V^j.
\tag7.1
$$
Then, using $Q_1,\,\ldots,\,Q_n$, one can write differential equations
\thetag{2.3} as follows:
$$
\hskip -2em
\left\{\aligned
&\dot x^i=V^i(x^1,\ldots,x^n,p_1,\ldots,p_n),\\
&\nabla_{\!t}p_i=Q_i(x^1,\ldots,x^n,p_1,\ldots,p_n).
\endaligned\right.
\tag7.2
$$
Components of covector $\bold Q$ in \thetag{7.2} play the same role as
components of vector $\bold F$ in \thetag{1.4}. Therefore covector field
$\bold Q$ is called {\bf force covector} or {\bf force field} of Newtonian
dynamical system \thetag{7.2}.
\head
8. Weak normality condition and weak normality equations.
\endhead
    Let's consider one-parametric family of trajectories of Newtonian
dynamical system \thetag{7.2}. This is a family of parametric curves
$q=q(t,y)$ in $T^*\!M$, where $t$ is time variable and $y$ is additional
parameter. In local chart this family of curves is represented by the
following set of $2n$ functions
$$
\xalignat 2
&\hskip -2em
\cases
x^1=x^1(t,y),\\
.\ .\ .\ .\ .\ .\ .\ .\ .\ .\\
x^n=x^n(t,y),
\endcases
&&\cases
p_1=p_1(t,y),\\
.\ .\ .\ .\ .\ .\ .\ .\ .\ .\\
p_n=p_n(t,y),
\endcases
\tag8.1
\endxalignat
$$
Differentiating first part of these functions \thetag{8.1} with
respect to additional parameter $y$, we obtain vector-function
$\boldsymbol\tau=\boldsymbol\tau(t,y)$ with components
$$
\hskip -2em
\tau^i=\frac{\partial x^i}{\partial y}.
\tag8.2
$$
Functions $p_1,\,\ldots,\,p_n$ in \thetag{8.1} are components of
covector-function. Therefore we should apply covariant derivative
to them:
$$
\hskip -2em
\xi_i=\nabla_{\!\boldsymbol\tau}p_i=\frac{\partial p_i}{\partial y}
-\sum^n_{j=1}\sum^n_{k=1}\Gamma^k_{ij}\,p_k\,\tau^j.
\tag8.3
$$
As a result we get covector-function $\boldsymbol\xi=\boldsymbol
\xi(t,y)$ with components \thetag{8.3}. Vector $\boldsymbol\tau$
is called {\bf variation vector} or, more exactly, {\bf vector of
variation of trajectory}, while covector $\boldsymbol\xi$ is called
{\bf covector of variation of momentum}.\par
    Components of both functions $\boldsymbol\tau(t,y)$ and
$\boldsymbol\xi(t,y)$ satisfy a system of ordinary differential
equations with respect to time variable. In order to derive these
equations let's differentiate equations \thetag{7.2} with respect to
parameter $y$. This yields
$$
\allowdisplaybreaks
\align
&\hskip -2em
\nabla_{\!t}\tau^i=\sum^n_{k=1}\nabla_{\!k}V^i\cdot\tau^k+
\sum^n_{k=1}\tilde\nabla^kV^i\cdot\xi_k,
\tag8.4\\
\vspace{2ex}
&\hskip -2em
\aligned
&\nabla_{\!t}\xi_i+\sum^n_{k=1}\!\left(\,\shave{\sum^n_{j=1}}
\shave{\sum^n_{s=1}}R^s_{ijk}\,p_s\,V^j-\shave{\sum^n_{j=1}}
\shave{\sum^n_{s=1}}D^{sj}_{ik}\,p_s\,Q_j\right)\cdot\tau^k\,+\\
&+\,\sum^n_{k=1}\sum^n_{j=1}\sum^n_{s=1}D^{sk}_{ij}\,p_s\,V^j\cdot
\xi_k=\sum^n_{k=1}\nabla_{\!k}Q_i\cdot\tau^k+
\sum^n_{k=1}\tilde\nabla^kQ_i\cdot\xi_k.
\endaligned
\tag8.5
\endalign
$$
Here in \thetag{8.5} we used components of curvature tensors
$\bold R$ and $\bold D$. First is an analog of standard curvature tensor
of Riemannian geometry, it is given by formula
$$
\hskip -2em
\gathered
R^k_{rij}=\frac{\partial\Gamma^k_{jr}}{\partial x^i}-
\frac{\partial\Gamma^k_{ir}}{\partial x^j}+\sum^n_{m=1}
\Gamma^k_{im}\,\Gamma^m_{jr}-\sum^n_{m=1}\Gamma^k_{jm}\,
\Gamma^m_{ir}\,+\\
\vspace{1ex}
+\sum^n_{m=1}\sum^n_{\alpha=1}p_\alpha\,\Gamma^\alpha_{mi}\,
\frac{\partial\Gamma^k_{jr}}{\partial p^m}-\sum^n_{m=1}
\sum^n_{\alpha=1}p_\alpha\,\Gamma^\alpha_{mj}\,\frac{\partial
\Gamma^k_{ir}}{\partial p^m}.
\endgathered
\tag8.6
$$
Second is a tensor of dynamic curvature, it is nonzero only for
extended connections, when $\Gamma^k_{ij}$ do actually depend on
components of momentum covector:
$$
\hskip -2em
D^{kr}_{ij}=-\frac{\partial\Gamma^k_{ij}}{\partial p_r}.
\tag8.7
$$\par
     In the next step we consider scalar product of vector
$\boldsymbol\tau$ and momentum covector $\bold p$. This scalar
product introduced by formula \thetag{1.13} determines so called
{\bf deviation function} $\varphi$. Like $\boldsymbol\tau=\boldsymbol
\tau(t,y)$ and $\boldsymbol\xi=\boldsymbol\xi(t,y)$, this is a function
of time variable $t$ and additional parameter $y$. In general case
deviation function \thetag{1.13} satisfies linear homogeneous ODE of
the order $2n$ (see theorem~6.1 in \cite{20}). But here we consider
special case determined by the following definition.
\definition{Definition 8.1} We say that Newtonian dynamical system
\thetag{2.3} satisfies {\bf weak} normality condition if for each
its trajectory $q=q(t)$ there is some second order homogeneous
linear ordinary differential equation 
$$
\hskip -2em
\ddot\varphi=\Cal A(t)\,\dot\varphi+\Cal B(t)\,\varphi
\tag8.8
$$
such that all deviation functions on the trajectory satisfy this
differential equation.
\enddefinition
   Saying {\it ``all deviation functions''}, in definition~8.1 we
imply that each trajectory $q=q(t)$ can be included into one-parametric
family of trajectories by various possible ways. Each such inclusion
defines some variation vector $\boldsymbol\tau=\boldsymbol\tau(t)$
and corresponding deviation function $\varphi=\varphi(t)$ on that
trajectory. Functions $\Cal A(t)$ and $\Cal B(t)$ in \thetag{8.8}
depend on the trajectory $q=q(t)$, but they do not depend on how
this trajectory is included into one parametric family of trajectories.
\par
    Let's calculate time derivatives of deviation function \thetag{1.13}.
For first order time derivative $\dot\varphi$ we obtain the following
expression:
$$
\hskip -2em
\dot\varphi=\nabla_{\!t}\varphi=\sum^n_{i=1}\nabla_{\!t}\tau^i\,p_i+
\sum^n_{i=1}\tau^i\,\nabla_{\!t}p_i.
\tag8.9
$$
Then we substitute \thetag{8.4} for $\nabla_{\!t}\tau^i$ and \thetag{7.2}
for $\nabla_{\!t}p_i$ into \thetag{8.9}. \pagebreak This yields more
detailed expression for first order time derivative of deviation function:
$$
\hskip -2em
\dot\varphi=\sum^n_{k=1}\left(\,\shave{\sum^n_{i=1}}
\nabla_{\!k}V^i\,p_i+Q_k\!\right)\cdot\tau^k+\sum^n_{k=1}
\left(\,\shave{\sum^n_{i=1}}\tilde\nabla^kV^i\,p_i\right)
\cdot\xi_k.
\tag8.10
$$
Applying $\nabla_{\!t}$ to \thetag{8.10}, we derive formula for second
order time derivative $\ddot\varphi$:
$$
\gather
\ddot\varphi=\sum^n_{k=1}\left(\,\shave{\sum^n_{i=1}}
\nabla_{\!k}V^i\,p_i+Q_k\!\right)\cdot\nabla_{\!t}\tau^k
+\sum^n_{k=1}\left(\,\shave{\sum^n_{i=1}}\tilde\nabla^kV^i
\,p_i\right)\cdot\nabla_{\!t}\xi_k\,+\\
+\,\sum^n_{k=1}\sum^n_{r=1}\left(\,\shave{\sum^n_{i=1}}
\nabla_{\!r}\!\nabla_{\!k}V^i\,p_i+\nabla_{\!r}Q_k\!\right)V^r\cdot
\tau^k+\sum^n_{k=1}\left(\,\shave{\sum^n_{i=1}}
\nabla_{\!k}V^i\,Q_i\right)\cdot\tau^k\,+\\
+\,\sum^n_{k=1}\sum^n_{r=1}\left(\,\shave{\sum^n_{i=1}}
\tilde\nabla^r\nabla_{\!k}V^i\,p_i+\tilde\nabla^r\!Q_k\!\right)Q_r
\cdot\tau^k
+\sum^n_{k=1}\left(\,\shave{\sum^n_{i=1}}\tilde\nabla^kV^i\,Q_i\right)
\cdot\xi_k\,+\\
+\,\sum^n_{k=1}\sum^n_{r=1}\left(\,\shave{\sum^n_{i=1}}
\nabla_{\!r}\!\tilde\nabla^kV^i\,p_i\right)V^r\cdot\xi_k
+\sum^n_{k=1}\sum^n_{r=1}\left(\,\shave{\sum^n_{i=1}}
\tilde\nabla^r\tilde\nabla^kV^i\,p_i\right)Q_r\cdot\xi_k.
\endgather
$$
Let's substitute \thetag{8.4} for $\nabla_{\!t}\tau^k$ and \thetag{8.5}
for $\nabla_{\!t}\xi_k$ into the above equality. It's easy to note that
resulting expression for $\ddot\varphi$ will be of the form
$$
\hskip -2em
\ddot\varphi=\sum^n_{k=1}\alpha^k\,\xi_k+\sum^n_{k=1}\beta_k\,\tau^k.
\tag8.11
$$
Here $\alpha^1,\,\ldots,\,\alpha^n$ are components of extended vector
field $\boldsymbol\alpha$ given by formula
$$
\gathered
\alpha^k=\sum^n_{i=1}\tilde\nabla^kV^i\,Q_i
+\sum^n_{r=1}\sum^n_{i=1}\nabla_{\!r}\!\tilde\nabla^kV^i\,p_i\,V^r
+\sum^n_{r=1}\sum^n_{i=1}\tilde\nabla^r\tilde\nabla^kV^i\,p_i\,Q_r\,+\\
+\sum^n_{r=1}\sum^n_{i=1}\tilde\nabla^rV^i\,p_i\,\tilde\nabla^kQ_r
-\sum^n_{r=1}\sum^n_{s=1}\sum^n_{i=1}\sum^n_{j=1}D^{sk}_{rj}\,
\tilde\nabla^rV^i\,p_i\,p_s\,V^j\,+\\
+\sum^n_{r=1}\sum^n_{i=1}\nabla_{\!r}V^i\,p_i\tilde\nabla^kV^r
+\sum^n_{r=1}\tilde\nabla^kV^r\,Q_r.
\endgathered\quad
\tag8.12
$$
Similarly, quantities $\beta_1,\,\ldots,\,\beta_n$ are components
of extended covector field $\boldsymbol\beta$:
$$
\hskip -2em
\gathered
\beta_k=\sum^n_{r=1}\sum^n_{i=1}\nabla_{\!r}\!\nabla_{\!k}V^i\,p_i
\,V^r+\sum^n_{r=1}V^r\,\nabla_{\!r}Q_k+\sum^n_{i=1}\nabla_{\!k}V^i
\,Q_i\,+\\
+\sum^n_{r=1}\sum^n_{i=1}\tilde\nabla^r\nabla_{\!k}V^i\,p_i\,Q_r
+\sum^n_{r=1}\tilde\nabla^r\!Q_k\,Q_r
+\sum^n_{r=1}\nabla_{\!k}V^r\,Q_r\,+\\
+\sum^n_{r=1}\sum^n_{i=1}\nabla_{\!r}V^i\,p_i\,\nabla_{\!k}V^r
+\sum^n_{r=1}\sum^n_{i=1}\tilde\nabla^rV^i\,p_i\,\nabla_{\!k}Q_r\,-\\
-\sum^n_{r=1}\sum^n_{s=1}\sum^n_{i=1}\sum^n_{j=1}\left(R^s_{rjk}\,
\tilde\nabla^rV^i\,p_i\,p_s\,V^j
-D^{sj}_{rk}\,\tilde\nabla^rV^i\,p_i\,p_s\,Q_j\right)
\endgathered
\tag8.13
$$
Let's compare formula \thetag{8.10} for time derivative $\dot\varphi$
with our previous notations \thetag{6.1}. It's easy to see that
\thetag{8.10} can be written as
$$
\hskip -2em
\dot\varphi=\sum^n_{k=1}\left(\,\shave{\sum^n_{i=1}}
\nabla_{\!k}V^i\,p_i+Q_k\!\right)\cdot\tau^k+\sum^n_{k=1}
W^k\cdot\xi_k.
\tag8.14
$$
Second term in \thetag{8.14} is scalar product of vector $\bold W$
and covector $\boldsymbol\xi$. Suppose that regularity condition
is fulfilled (see definition~6.1). Then we can introduce same scalar
product into the formula \thetag{8.11} for second order partial
derivative $\ddot\varphi$:
$$
\hskip -2em
\ddot\varphi=\sum^n_{k=1}\frac{\left<\bold p\,|\,\boldsymbol
\alpha\right>}{\Omega}\,W^k\cdot\xi_k+\sum^n_{k=1}\sum^n_{r=1}
\alpha^r\,P^k_r\cdot\xi_k+\sum^n_{k=1}\beta_k\,\tau^k.
\tag8.15
$$
Here we used formula \thetag{6.3} written in the following form:
$$
\hskip -2em
\delta^k_r=P^k_r+\frac{W^k\,p_r}{\Omega}.
\tag8.16
$$
Combining formulas \thetag{8.14} and \thetag{8.15} for time derivatives
$\dot\varphi$ and $\ddot\varphi$, we obtain the equality with right hand
side free of scalar product $\left<\boldsymbol\xi\,|\,\bold W\right>$:
$$
\hskip -2em
\gathered
\ddot\varphi-\frac{\left<\bold p\,|\,\boldsymbol\alpha\right>}
{\Omega}\,\dot\varphi=\sum^n_{k=1}\sum^n_{r=1}\alpha^r\,P^k_r
\cdot\xi_k+\sum^n_{k=1}\beta_k\cdot\tau^k\,-\\
-\,\sum^n_{k=1}\left(\,\shave{\sum^n_{i=1}}\shave{\sum^n_{s=1}}
\frac{\nabla_{\!k}V^i\,p_i\,\alpha^s\,p_s}{\Omega}
+\shave{\sum^n_{s=1}}\frac{Q_k\,\alpha^s\,p_s}{\Omega}\right)
\cdot\tau^k.
\endgathered
\tag8.17
$$
For the sake of brevity in further calculations we introduce the
following notations:
$$
\hskip -2em
\eta_k=\beta_k-\sum^n_{i=1}\sum^n_{s=1}\frac{\nabla_{\!k}V^i\,p_i\,
\alpha^s\,p_s}{\Omega}-\sum^n_{s=1}\frac{Q_k\,\alpha^s\,p_s}{\Omega}.
\tag8.18
$$
Quantities $\eta_1,\,\ldots,\,\eta_n$ in \thetag{8.18} are components
of covector field $\boldsymbol\eta$. Using notations \thetag{8.18}, we
can simplify the above equality \thetag{8.17}:
$$
\hskip -2em
\ddot\varphi-\frac{\left<\bold p\,|\,\boldsymbol\alpha\right>}
{\Omega}\,\dot\varphi=\sum^n_{k=1}\sum^n_{r=1}\alpha^r\,P^k_r
\cdot\xi_k+\sum^n_{k=1}\eta_k\cdot\tau^k.
\tag8.19
$$
Deviation function $\varphi$ is scalar product of vector $\boldsymbol
\tau$ and covector $\bold p$. Using \thetag{8.16}, we can introduce
such scalar product into right hand side of \thetag{8.19}. Then we get
$$
\ddot\varphi-\frac{\left<\bold p\,|\,\boldsymbol\alpha\right>}
{\Omega}\,\dot\varphi-\frac{\left<\boldsymbol\eta\,|\,\bold W\right>}
{\Omega}\,\varphi=\sum^n_{k=1}\sum^n_{r=1}\alpha^r\,P^k_r
\cdot\xi_k+\sum^n_{k=1}\sum^n_{r=1}\eta_r\,P^r_k\cdot\tau^k.
\quad
\tag8.20
$$
Now let's recall the definition~8.1. Suppose that the following
equalities are fulfilled:
$$
\xalignat 2
\hskip -2em
&\sum^n_{r=1}\alpha^r\,P^k_r=0,
&&\sum^n_{r=1}\eta_r\,P^r_k=0.
\tag8.21
\endxalignat
$$
Then deviation function $\varphi$ satisfies second order ordinary
differential equation \thetag{8.8} with coefficients $\Cal A$ and
$\Cal B$ determined by vector field $\boldsymbol\alpha$ and
covector field $\boldsymbol\eta$:
$$
\xalignat 2
\hskip -2em
&\Cal A=\frac{\left<\bold p\,|\,\boldsymbol\alpha\right>}{\Omega},
&&\Cal B=\frac{\left<\boldsymbol\eta\,|\,\bold W\right>}
{\Omega}.
\tag8.22
\endxalignat
$$
This means that equations \thetag{8.21} are sufficient for weak
normality condition formulated in definition~8.1 to be fulfilled.
\par
    Conversely, suppose that weak normality condition for Newtonian
dynamical system \thetag{7.2} is fulfilled. Combining \thetag{8.8}
with \thetag{8.20}, we derive
$$
\hskip -2em
\tilde\Cal A\,\dot\varphi+\tilde\Cal B\,\varphi
=\sum^n_{k=1}\sum^n_{r=1}\alpha^r\,P^k_r\cdot\xi_k
+\sum^n_{k=1}\sum^n_{r=1}\eta_r\,P^r_k\cdot\tau^k,
\tag8.23
$$
\vskip -2ex
\noindent where
\vskip -4ex
$$
\xalignat 2
&\hskip -2em
\tilde\Cal A=\Cal A-\frac{\left<\bold p\,|\,\boldsymbol\alpha\right>}
{\Omega},
&&\tilde\Cal B=\Cal B-\frac{\left<\boldsymbol\eta\,|\,\bold W\right>}
{\Omega}.
\tag8.24
\endxalignat
$$
Substituting \thetag{8.14} and \thetag{1.13} for $\dot\varphi$ and
$\varphi$ into \thetag{8.23}, we obtain
$$
\hskip -2em
\tilde\Cal A\,\left<\boldsymbol\xi\,|\,\bold W\right>+
\tilde\Cal C\,\left<\bold p\,|\,\boldsymbol\tau\right>
=\sum^n_{k=1}\sum^n_{r=1}\alpha^r\,P^k_r\cdot\xi_k
+\sum^n_{k=1}\sum^n_{r=1}\tilde\eta_r\,P^r_k\cdot\tau^k.
\tag8.25
$$
Here for more convenience we introduced the following auxiliary
notations:
$$
\align
&\hskip -2em
\tilde\Cal C=\tilde\Cal B+\tilde\Cal A\left(\,\shave{\sum^n_{r=1}}
\shave{\sum^n_{i=1}}\frac{\nabla_{\!r}V^i\,p_i\,W^r}{\Omega}
+\shave{\sum^n_{r=1}}\frac{W^r\,Q_r}{\Omega}\!\right),
\tag8.26\\
&\hskip -2em
\tilde\eta_r=\eta_r-\tilde\Cal A\left(\,\shave{\sum^n_{i=1}}
\nabla_{\!r}V^i\,p_i+Q_r\!\right).
\tag8.27
\endalign
$$\par 
     Now conceptual point\,! Looking at formula \thetag{8.25},
let's remember that for a fixed trajectory $p=p(t)$ of Newtonian
dynamical system \thetag{7.2} functions $\tau^1(t),\,\ldots,\,
\tau^n(t)$ and $\xi_1(t),\,\ldots,\,\xi_n(t)$ satisfy system of
first order linear homogeneous ordinary differential equations
(see \thetag{8.4} and \thetag{8.5}). Solutions of these differential
equations form $2n$-dimensional linear space $\goth T$, while
components of vector-function $\boldsymbol\tau(t)$ and
covector-function $\boldsymbol \xi(t)$ for any fixed time instant
$t=t_0$ can be treated as coordinates in $\goth T$. In other words,
linear $\goth T$ is isomorphic to direct sum
$$
\hskip -2em
\goth T\cong T_p(M)\oplus T^*_p(M),
\tag8.28
$$
where $p=p(t_0)$. In regular case (see definition~6.1) denominator
$\Omega=\Omega(p,\bold p)$ defined by formula \thetag{6.2} is nonzero
for $\bold p\neq 0$. Hence vector $\bold W=\bold W(p,\bold p)$ is
transversal to null-space of momentum covector $\bold p$:
$$
\hskip -2em
\goth P=\Img\bold P=\left\{\bold X\in T_p(M):\ \left<\bold p\,|\,\bold X
\right>=0\right\}.
\tag8.29
$$
Therefore tangent space $T_p(P)$ in \thetag{8.28} breaks into direct sum
$$
\hskip -2em
T_p(M)=\left<\bold W\right>\oplus \goth P,
\tag8.30
$$
where $\left<\bold W\right>$ is linear span of vector $\bold W$. We can
write similar expansion for $T^*_p(M)$:
$$
\hskip -2em
T^*_p(M)=\left<\bold p\right>\oplus \goth W.
\tag8.31
$$
Here $\left<\bold p\right>$ is linear span of covector $\bold p$ and
$\goth W$ is null-space of $\bold W$:
$$
\hskip -2em
\goth W=\Img\bold P^*=\left\{\bold y\in T^*_p(M):\ \left<\bold y\,|\,
\bold W\right>=0\right\}.
\tag8.32
$$
Substituting \thetag{8.30} and \thetag{8.31} into \thetag{8.28} we
get the following expansion:
$$
\hskip -2em
\goth T\cong\underbrace{\left<\bold W\right>\oplus
\left<\bold p\right>}\oplus\underbrace{\goth P\oplus\goth W}=
\goth T_1\oplus\goth T_2.
\tag8.33
$$
Expansion \thetag{8.33} of the space $\goth T$ generates conjugate
expansion of dual space $\goth T^*$:
$$
\hskip -2em
\goth T^*\cong\goth T^*_1\oplus\goth T^*_2.
\tag8.34
$$
The above equality \thetag{8.25} is an equality of functions. However,
if we treat pair of $\boldsymbol\tau$ and $\boldsymbol\xi$ as a point
of $\goth T$, then for any fixed $t=t_0$ both sides of the equality
\thetag{8.25} can be treated as linear functionals in $\goth T$, i\.\,
e\. they are elements of $\goth T^*$. Note that left hand side of
\thetag{8.25} is in $\goth T^*_1$, while right hand side is an element
of $\goth T^*_2$. Subspaces $\goth T^*_1$ and $\goth T^*_2$ in direct
sum \thetag{8.34} have zero intersection. Therefore both sides of
\thetag{8.25} are zero (as elements of $\goth T^*$ and as functions
of $t$ as well, since $t=t_0$ as an arbitrary fixed instant of time).
Thus we have the equality
$$
\hskip -2em
\tilde\Cal A\,\left<\boldsymbol\xi\,|\,\bold W\right>+
\tilde\Cal C\,\left<\bold p\,|\,\boldsymbol\tau\right>=0.
\tag8.35
$$
Due to expansion $\goth T^*_1\cong\left<\bold W\right>^*\oplus
\left<\bold p\right>^*$ from the equality \thetag{8.35} we derive
$$
\xalignat 2
&\hskip -2em
\tilde\Cal A=0,
&&\tilde\Cal C=0.
\tag8.36
\endxalignat
$$
Applying \thetag{8.26}, \thetag{8.27}, and \thetag{8.24} to
\thetag{8.36}, we get $\tilde\Cal B=0$, $\tilde\eta_k=\eta_k$
and come back to \thetag{8.22}. Substituting $\tilde\Cal A=0$
and $\tilde\Cal B=0$ into \thetag{8.23}, we obtain the equality
$$
\sum^n_{k=1}\sum^n_{r=1}\alpha^r\,P^k_r\cdot\xi_k
+\sum^n_{k=1}\sum^n_{r=1}\eta_r\,P^r_k\cdot\tau^k=0,
$$
which breaks  into two separate parts equivalent to \thetag{8.21}.
Thus, assuming that weak normality condition is fulfilled, we have
derived again the equalities \thetag{8.21}. When written explicitly,
they form a system of partial differential equations relating
extended vector field $\bold V$, extended covector field $\bold Q$,
and extended affine connection $\Gamma$. These equations are called
{\bf weak normality equations}. Just above we have proved the following
theorem for them.
\proclaim{Theorem 8.1} Weak normality condition for Newtonian
dynamical system \thetag{7.2} is equivalent to weak normality
equations \thetag{8.21} that should be fulfilled at all
points $q=(p,\bold p)$ of cotangent bundle $T^*\!M$, where
$\bold p\neq 0$.
\endproclaim
\noindent Using \thetag{8.12}, \thetag{8.13}, \thetag{8.18}, one can
easily write weak normality equations \thetag{8.21} explicitly. However,
we shall not do it, since they are rather huge.
\head
9. Additional normality condition.
\endhead
   Let's proceed with studying normal shift phenomenon assuming that
weak normality condition for Newtonian dynamical system \thetag{7.2}
is fulfilled. For this purpose let's consider some hypersurface
$\sigma$ and apply initial data \thetag{1.14} to differential
equations \thetag{7.2}. As a result we get $(n-1)$-parametric family
of trajectories of Newtonian dynamical system \thetag{7.2}. Let's
fix some point $p_0\in \sigma$. If $y^1,\,\ldots,y^{n-1}$ are
local coordinates on $\sigma$ in some neighborhood of the point $p_0$
and if $x^1,\,\ldots,x^n$ are local coordinates in $M$ in $n$-dimensional
neighborhood of this point, then our $(n-1)$-parametric family of
trajectories is represented by functions
$$
\xalignat 2
&\hskip -2em
\cases
x^1=x^1(t,y^1,\ldots,y^{n-1}),\\
.\ .\ .\ .\ .\ .\ .\ .\ .\ .\ .\ .\ .\ .\ .\ .\ .\ .\\
x^n=x^n(t,y^1,\ldots,y^{n-1}),
\endcases
&&\cases
p_1=p_1(t,y^1,\ldots,y^{n-1}),\\
.\ .\ .\ .\ .\ .\ .\ .\ .\ .\ .\ .\ .\ .\ .\ .\ .\ .\\
p_n=p_n(t,y^1,\ldots,y^{n-1}).
\endcases
\tag9.1
\endxalignat
$$
Comparing \thetag{9.1} with \thetag{8.1} and looking at \thetag{8.2}
and \thetag{8.3}, we see that now we can define vector-functions
$\boldsymbol\tau_1,\,\ldots,\,\boldsymbol\tau_{n-1}$ and
covector-functions $\boldsymbol\xi_1,\,\ldots,\,\boldsymbol\xi_{n-1}$.
Their components in local chart are given by the following derivatives:
$$
\xalignat 2
&\hskip -2em
\tau^s_i=\frac{\partial x^s}{\partial y^i},
&&\xi_{si}=\nabla_{\!\boldsymbol\tau_i}p_s.
\tag9.2
\endxalignat
$$
Vector-functions $\boldsymbol\tau_1,\,\ldots,\,\boldsymbol\tau_{n-1}$
with components \thetag{9.2} define $n-1$ deviation functions
$\varphi_1,\,\ldots,\,\varphi_{n-1}$ according to the formula
\thetag{1.13}:
$$
\hskip -2em
\varphi_i=\varphi_i(t,y^1,\ldots,y^{n-1})=\left<\bold p\,|\,
\boldsymbol\tau_i\right>=\sum^n_{s=1}\tau^s_i\,p_s.
\tag9.3
$$
We assumed that weak normality condition is fulfilled (see definition~8.1).
Therefore each function $\varphi_i$ in \thetag{9.3} satisfies differential
equation of the form \thetag{8.8}:
$$
\hskip -2em
\ddot\varphi_i=\Cal A(t,y^1,\ldots,y^{n-1})\,\dot\varphi_i
+\Cal B(t,y^1,\ldots,y^{n-1})\,\varphi_i.
\tag9.4
$$
According to definition~1.4, in order to have normal shift we should
provide vanishing of all deviation functions $\varphi_i,\,\ldots,\,
\varphi_{n-1}$ in \thetag{9.3}. Due to the equation \thetag{9.4} it
is sufficient to provide the following initial data for them:
$$
\xalignat 2
&\hskip -2em
\varphi_i\,\hbox{\vrule height 8pt depth 8pt width 0.5pt}_{\,t=0}=0,
&&\dot\varphi_i\,\hbox{\vrule height 8pt depth 8pt width 0.5pt}_{\,t=0}=0.
\tag9.5
\endxalignat
$$
First part of initial conditions \thetag{9.5} is fulfilled due to
initial data \thetag{1.14}. Second part of these initial conditions
should be fulfilled at the expense of proper choice of scalar function
$\nu=\nu(p)=\nu(y^1,\ldots,y^n)$ in \thetag{1.14}. In order to calculate
time derivative $\dot\varphi_i$ we can use formula \thetag{8.14}. Here
it is written as follows:
$$
\hskip -2em
\dot\varphi_i=\sum^n_{s=1}\left(\,\shave{\sum^n_{r=1}}
\nabla_{\!s}V^r\,p_r+Q_s\!\right)\cdot\tau^s_i+\sum^n_{s=1}
W^s\cdot\xi_{si}.
\tag9.6
$$
Then we should substitute $t=0$ into \thetag{9.6}. Vectors $\boldsymbol
\tau_1,\,\ldots,\,\boldsymbol\tau_{n-1}$ at initial instant of time
$t=0$ form base of coordinate vectors in tangent space to initial
hypersurface $\sigma$. Initial value of momentum covector $\bold p$
is given by \thetag{1.14}. Thus we have to calculate initial values
for covectors $\boldsymbol\xi_1,\,\ldots,\,\boldsymbol\xi_{n-1}$
in formula \thetag{9.6}:
$$
\hskip -2em
\xi_{si}\,\hbox{\vrule height 8pt depth 8pt width 0.5pt}_{\,t=0}
=\nabla_{\!\boldsymbol\tau_i}p_s\,\hbox{\vrule height 8pt depth
8pt width 0.5pt}_{\,t=0}=\frac{\partial\nu}{\partial y^i}\cdot
n_s+\nu\cdot\nabla_{\!\boldsymbol\tau_i}n_s.
\tag9.7
$$
Substituting \thetag{9.7} into \thetag{9.6} and then substituting
\thetag{9.6} into \thetag{9.5}, we get
$$
\hskip -2em
\frac{\partial\nu}{\partial y^i}=-\frac{(\nu)^2}{\Omega}
\sum^n_{s=1}W^s\cdot\nabla_{\!\boldsymbol\tau_i}n_s-\frac{\nu}
{\Omega}\sum^n_{s=1}\left(\,\shave{\sum^n_{r=1}}
\nabla_{\!s}V^r\,p_r+Q_s\!\right)\cdot\tau^s_i.
\tag9.8
$$
In two-dimensional case $n=\dim M=2$ we have only one variable $y=y^1$.
Then \thetag{9.8} turns to ordinary differential equation with respect
to function $\nu=\nu(y)$. Normalizing condition \thetag{1.3} yields
initial value problem for this ordinary differential equation, which is
always solvable (at least, locally).\par
    In multidimensional case $n\geqslant 3$ we have several variables
$y^1,\,\ldots,\,y^{n-1}$. In this case partial differential equations
\thetag{9.8} form complete system of Pfaff equations with respect to
function $\nu=\nu(y^1,\,\ldots,\,y^{n-1})$. Initial value problem
$$
\hskip -2em
\nu(p_0)=\nu_0
\tag9.9
$$
at some fixed point $p_0\in\sigma$ with local coordinates $y^1_0,
\,\ldots,\,y^{n-1}_0$ is typical for Pfaff equations \thetag{9.8}.
However, now it is not unconditionally solvable (even locally).
\definition{Definition 9.1} Complete system of Pfaff equations
\thetag{9.8} is called {\it compatible} if initial value problem
\thetag{9.9} for it is locally solvable for all $p_0\in M$ and for
all $\nu_0\neq 0$.
\enddefinition
\noindent Let's write Pfaff equations \thetag{9.8} formally,
denoting by $\psi_i$ their right hand sides:
$$
\hskip -2em
\frac{\partial\nu}{\partial y^i}=\psi_i(\nu,y^1,\ldots,y^{n-1})
\tag9.10
$$
Due to \thetag{9.10} we can calculate mixed partial derivatives
of $\nu$ in two different ways
$$
\align
\hskip -2em\frac{\partial^2\nu}{\partial y^i\,\partial y^j}&=
\frac{\partial\psi_i}{\partial y^j}+\frac{\partial\psi_i}{\partial\nu}
\,\psi_j=\vartheta_{ij}(\nu,y^1,\ldots,y^{n-1}),\hskip -2em
\tag9.11\\
\vspace{1ex}
\hskip -2em\frac{\partial^2\nu}{\partial y^j\,\partial y^i}&=
\frac{\partial\psi_j}{\partial y^i}+\frac{\partial\psi_j}{\partial\nu}
\,\psi_i=\vartheta_{ji}(\nu,y^1,\ldots,y^{n-1}).\hskip -2em
\tag9.12
\endalign
$$
Equating \thetag{9.11} and \thetag{9.12}, we get compatibility condition
for \thetag{9.10}:
$$
\hskip -2em
\vartheta_{ij}(\nu,y^1,\ldots,y^{n-1})=\vartheta_{ji}(\nu,y^1,
\ldots,y^{n-1}).
\tag9.13
$$
\proclaim{Lemma 9.1} Pfaff equations \thetag{9.10} are compatible in
the sense of definition~9.1 if and only if for $\nu\neq 0$ left and
right hands sides of \thetag{9.13} are equal to each other identically
as functions of $n$ independent variables $y^1,\,\ldots,\,y^{n-1}$,
and $\nu$.
\endproclaim
Lemma~9.1 is standard result in the theory of Pfaff equations. Proof
of this lemma can be found in thesis \cite{7}.
\definition{Definition 9.2} We say that Newtonian dynamical system
satisfies {\bf additional} normality condition if Pfaff equations
for the function $\nu(p)$ in \thetag{1.14} derived from initial
conditions \thetag{9.5} are compatible for any hypersurface $\sigma$
in $M$.
\enddefinition
\noindent For the sake of convenience we introduce extended covector
field $\bold U$ with components
$$
\hskip -2em
U_s=\sum^n_{r=1}\nabla_{\!s}V^r\,p_r+Q_s.
\tag9.14
$$
Then Pfaff equations \thetag{9.8} are written as follows:
$$
\hskip -2em
\frac{\partial\nu}{\partial y^i}=-\frac{(\nu)^2}{\Omega}
\sum^n_{s=1}W^s\cdot\nabla_{\!\boldsymbol\tau_i}n_s-\frac{\nu}
{\Omega}\sum^n_{s=1}U_s\cdot\tau^s_i.
\tag9.15
$$
Now, relying upon definition~9.2, we shall derive explicit form of
compatibility equation \thetag{9.13}. For this purpose let's calculate
partial derivatives \thetag{9.11} and \thetag{9.12} in explicit form.
For partial derivatives \thetag{9.11} we obtain
$$
\gathered
\frac{\partial^2\nu}{\partial y^i\,\partial y^j}=\nabla_{\!\boldsymbol
\tau_j}\psi_i=\frac{2\,(\nu)^3}{\Omega^2}\sum^n_{s=1}\sum^n_{r=1}W^s\,W^r
\,\nabla_{\!\boldsymbol\tau_i}n_s\,\nabla_{\!\boldsymbol\tau_j}n_r\,+\\
+\,\frac{2\,(\nu)^2}{\Omega^2}\sum^n_{s=1}\sum^n_{r=1}W^s\,U_r
\,\nabla_{\!\boldsymbol\tau_i}n_s\,\tau^r_j
+\frac{(\nu)^2}{\Omega^2}\sum^n_{s=1}\sum^n_{r=1}U_s\,W^r
\,\tau^s_i\,\nabla_{\!\boldsymbol\tau_j}n_r\,+\\
+\,\frac{\nu}{\Omega^2}\sum^n_{s=1}\sum^n_{r=1}U_s\,U_r\,\tau^s_i\,\tau^r_j
-(\nu)^2\sum^n_{s=1}\nabla_{\!\boldsymbol\tau_j}\!
\left(\frac{W^s}{\Omega}\right)\cdot\nabla_{\!\boldsymbol\tau_i}n_s\,-\\
-\,\nu\sum^n_{s=1}\nabla_{\!\boldsymbol\tau_j}\!\left(\frac{U_s}{\Omega}
\right)\,\tau^s_i
-\frac{(\nu)^2}{\Omega}\sum^n_{s=1}W^s\,
\nabla_{\!\boldsymbol\tau_j}\!\nabla_{\!\boldsymbol\tau_i}n_s
-\frac{\nu}{\Omega}\sum^n_{s=1}U_s\,\nabla_{\!\boldsymbol\tau_j}\tau^s_i.
\endgathered\quad
\tag9.16
$$
In order to transform \thetag{9.16} we need to bring about some
preliminary calculations. For covariant derivative $\nabla_{\!\boldsymbol
\tau_j}\tau^s_i$ in \thetag{9.16} we have
$$
\nabla_{\!\boldsymbol\tau_j}\tau^s_i=\frac{\partial\tau^s_i}{\partial
y^j}+\sum^n_{r=1}\sum^n_{q=1}\Gamma^s_{rq}\,\tau^r_i\,\tau^q_j=
\frac{\partial^2x^s}{\partial y^i\,\partial y^j}+\sum^n_{r=1}\sum^n_{q=1}
\Gamma^s_{rq}\,\tau^r_i\,\tau^q_j.
$$
Taking into account symmetry of connection components $\Gamma^s_{rq}
=\Gamma^s_{qr}$, we derive
$$
\hskip -2em
\nabla_{\!\boldsymbol\tau_i}\tau^s_j-\nabla_{\!\boldsymbol\tau_j}
\tau^s_i=0.
\tag9.17
$$
In a similar way by direct calculations we derive the following
identities:
$$
\hskip -2em
\aligned
&\nabla_{\!\boldsymbol\tau_j}\!\left(\frac{W^s}{\Omega}\right)
=\sum^n_{r=1}\nabla_{\!r}\!\left(\frac{W^s}{\Omega}\right)
\cdot\tau^r_j+\sum^n_{r=1}\tilde\nabla^r\!\left(\frac{W^s}{\Omega}\right)
\cdot\xi_{rj},\\
&\nabla_{\!\boldsymbol\tau_j}\!\left(\frac{U_s}{\Omega}\right)
=\sum^n_{r=1}\nabla_{\!r}\!\left(\frac{U_s}{\Omega}\right)
\cdot\tau^r_j+\sum^n_{r=1}\tilde\nabla^r\!\left(\frac{U_s}{\Omega}
\right)\cdot\xi_{rj}.
\endaligned
\tag9.18
$$
Formulas \thetag{9.18} are special cases of general formula applicable
to arbitrary extended tensor field $\bold X$. If $X^{i_1\ldots\,i_r}_{j_1
\ldots\,j_s}$ are components of $\bold X$ in local chart, then we have
$$
\nabla_{\!\boldsymbol\tau_q}X^{i_1\ldots\,i_r}_{j_1\ldots\,j_s}
=\sum^n_{r=1}\nabla_{\!r}X^{i_1\ldots\,i_r}_{j_1\ldots\,j_s}
\cdot\tau^r_q+\sum^n_{r=1}\tilde\nabla^rX^{i_1\ldots\,i_r}_{j_1
\ldots\,j_s}\cdot\xi_{rq}
$$
And finally, there is an identity for commutator of two covariant
derivatives:
$$
\hskip -2em
\gathered
[\nabla_{\!\boldsymbol\tau_i},\,\nabla_{\!\boldsymbol\tau_j}]n_s
=-\sum^n_{q=1}\sum^n_{\alpha=1}\sum^n_{\gamma=1}\frac{p_q}{\nu}\,
R^q_{s\alpha\gamma}\,\tau^\alpha_i\,\tau^\gamma_j\,+\\
+\sum^n_{r=1}\sum^n_{\alpha=1}
\sum^n_{\gamma=1}\frac{p_\alpha}{\nu}\,D^{\alpha\gamma}_{sr}\,
\tau^r_j\,\xi_{\gamma i}
-\sum^n_{r=1}\sum^n_{\alpha=1}
\sum^n_{\gamma=1}\frac{p_\alpha}{\nu}\,D^{\alpha\gamma}_{sr}\,
\tau^r_i\,\xi_{\gamma j}.
\endgathered
\tag9.19
$$
Combining \thetag{9.15} and \thetag{9.7}, for the quantities $\xi_{rj}$
in \thetag{9.18} we derive
$$
\hskip -2em
\xi_{rj}=\sum^n_{s=1}\nu\,P^s_r\,\nabla_{\!\boldsymbol\tau_j}n_s-
\sum^n_{s=1}\frac{U_s\,p_r}{\Omega}\,\tau^s_j.
\tag9.20
$$
Here $P^s_r$ are components of projection operator $\bold P$ introduced
in \thetag{6.3}. Further
$$
\tilde\nabla^r\!\left(\frac{W^s}{\Omega}\right)=\frac{1}{\Omega^2}\!
\left(\Omega\,\tilde\nabla^rW^s-W^s\,W^r-\shave{\sum^n_{q=1}}W^s\,
p_q\,\tilde\nabla^rW^q\right).
$$
Its is easy to see that right hand side of this formula simplifies 
when we introduce components of projector operator $\bold P$. Indeed,
we have
$$
\hskip -2em
\tilde\nabla^r\!\left(\frac{W^s}{\Omega}\right)=-\frac{W^s\,W^r}
{\Omega^2}+\sum^n_{q=1}\frac{\tilde\nabla^rW^q}{\Omega}\,P^s_q.
\tag9.21
$$
Combining \thetag{9.20} and \thetag{9.21}, we obtain the following
equality:
$$
\hskip -2em
\gathered
\sum^n_{r=1}\tilde\nabla^r\!\left(\frac{W^s}{\Omega}\right)\,\xi_{rj}=
\sum^n_{k=1}\sum^n_{r=1}\sum^n_{q=1}\nu\,P^k_r\,\frac{\tilde\nabla^rW^q}
{\Omega}\,P^s_q\,\nabla_{\!\boldsymbol\tau_j}n_k\,+\\
+\sum^n_{k=1}\frac{W^s}{\Omega^2}\,U_k\,\tau^k_j
-\sum^n_{k=1}\sum^n_{r=1}\sum^n_{q=1}\frac{p_r\,\tilde\nabla^rW^q}
{\Omega^2}\,P^s_q\,U_k\,\tau^k_j.
\endgathered
\tag9.22
$$
In similar way for the first term in right hand side of first
equality \thetag{9.18} we derive
$$
\hskip -2em
\sum^n_{r=1}\nabla_{\!r}\!\left(\frac{W^s}{\Omega}\right)
\cdot\tau^r_j=\sum^n_{r=1}\sum^n_{q=1}\frac{\nabla_{\!r}W^q}{\Omega}
\,P^s_q\,\tau^r_j.
\tag9.23
$$
Now we are able to proceed with transforming the equality \thetag{9.16}.
Note that terms symmetric in indices $i$ and $j$ make no contribution
to ultimate compatibility equation \thetag{9.13}. Therefore further we
shall omit them replacing by dots. Taking into account \thetag{9.17},
\thetag{9.18}, \thetag{9.19}, \thetag{9.22}, and \thetag{9.23}, for
\thetag{9.13} we derive
$$
\gathered
\theta_{ij}-\theta_{ji}=
(\nu)^3\sum^n_{k=1}\sum^n_{q=1}\sum^n_{r=1}\sum^n_{s=1}
\,\frac{\tilde\nabla^rW^s-\tilde\nabla^sW^r}
{\Omega}\,(P^q_r\,\nabla_{\!\boldsymbol\tau_i}n_q)\,
(P^k_s\,\nabla_{\!\boldsymbol\tau_j}n_k)\,+\\
+\,(\nu)^2\sum^n_{q=1}\sum^n_{r=1}\sum^n_{s=1}\left(
\frac{\tilde\nabla^rU_s}{\Omega}+\shave{\sum^n_{m=1}}
\frac{\tilde\nabla^mW^r-\tilde\nabla^rW^m}{\Omega^2}\,U_s\,p_m
-\frac{\nabla_{\!s}W^r}{\Omega}\,+\right.\\
\left.+\shave{\sum^n_{m=1}}\shave{\sum^n_{k=1}}\frac{W^k\,p_m
\,D^{mr}_{ks}}{\Omega}\right)(P^q_r\,\nabla_{\!\boldsymbol
\tau_i}n_q)\,\tau^s_j\,-(\nu)^2\sum^n_{q=1}\sum^n_{r=1}
\sum^n_{s=1}\left(\frac{\tilde\nabla^rU_s}{\Omega}\,+
\vphantom{\shave{\sum^n_{m=1}}}\right.
\\
\left.+\,\shave{\sum^n_{m=1}}\frac{\tilde\nabla^mW^r
-\tilde\nabla^rW^m}{\Omega^2}\,U_s\,p_m-\frac{\nabla_{\!s}W^r}
{\Omega}+\shave{\sum^n_{m=1}}\shave{\sum^n_{k=1}}\frac{W^k\,p_m
\,D^{mr}_{ks}}{\Omega}\right)\times\\
\times\,(P^q_r\,\nabla_{\!\boldsymbol\tau_j}n_q)\,\tau^s_i
+\nu\sum^n_{r=1}\sum^n_{s=1}\left(\,
\shave{\sum^n_{m=1}}
\frac{p_m\,\tilde\nabla^mU_r}{\Omega^2}\,U_s
-\shave{\sum^n_{m=1}}
\frac{p_m\,\tilde\nabla^mU_s}{\Omega^2}\,U_r\,+\right.
\\
\left.+\,\frac{\nabla_{\!r}U_s}{\Omega}-\frac{\nabla_{\!s}U_r}
{\Omega}+\shave{\sum^n_{m=1}}\frac{p_m\,\nabla_{\!s}W^m\,U_r}{\Omega^2}
-\shave{\sum^n_{m=1}}\frac{p_m\,\nabla_{\!r}W^m\,U_s}{\Omega^2}\right)
\tau^r_i\,\tau^s_j\,+
\\
+\,\nu\,\sum^n_{r=1}\sum^n_{s=1}\sum^n_{k=1}\sum^n_{q=1}W^k\,p_q
\left(\,\shave{\sum^n_{m=1}}\frac{D^{mq}_{kr}\,U_s-D^{mq}_{ks}\,U_r}
{\Omega^2}\,p_m-\frac{R^q_{krs}}{\Omega}\right)\tau^r_i\,\tau^s_j.
\endgathered
$$
In deriving the above equality we also used the following quite obvious
formulas:
$$
\hskip -2em
\aligned
&\tilde\nabla^r\!\left(\frac{U_s}{\Omega}\right)=
\frac{\tilde\nabla^rU_s}{\Omega}-\frac{W^r\,U_s}{\Omega^2}
-\sum^n_{m=1}\frac{p_m\,\tilde\nabla^rW^m\,U_s}{\Omega^2},\\
&\nabla_{\!r}\!\left(\frac{U_s}{\Omega}\right)=
\frac{\nabla_{\!r}U_s}{\Omega}
-\sum^n_{m=1}\frac{p_m\,\nabla_{\!r}W^m\,U_s}{\Omega^2}.
\endaligned
\tag9.24
$$
Using \thetag{9.24}, for two summands in right hand side of second
equality \thetag{9.18} we get
$$
\align
&\hskip -2em
\gathered
\sum^n_{r=1}\tilde\nabla^r\!\left(\frac{U_s}{\Omega}\right)\cdot
\xi_{rj}=-\nu\sum^n_{r=1}\sum^n_{k=1}\sum^n_{m=1}\frac{p_m\,\tilde
\nabla^rW^m\,U_s}{\Omega^2}\,P^k_r\,\nabla_{\!\boldsymbol\tau_j}
n_k\,+\\
+\,\nu\sum^n_{r=1}\sum^n_{k=1}\frac{\tilde\nabla^rU_s}{\Omega}
\,P^k_r\,\nabla_{\!\boldsymbol\tau_j}n_k-\sum^n_{k=1}\sum^n_{m=1}
\frac{p_m\,\tilde\nabla^mU_s}{\Omega^2}\,U_k\,\tau^k_j\,+\\
+\sum^n_{k=1}\frac{1}{\Omega^2}\,U_s\,U_k\,\tau^k_j+\sum^n_{k=1}
\sum^n_{q=1}\sum^n_{m=1}\frac{p_q\,p_m\,\tilde\nabla^qW^m}
{\Omega^3}\,U_s\,U_k\,\tau^k_j,
\endgathered
\tag9.25\\
&\hskip -2em
\gathered
\sum^n_{r=1}\nabla_{\!r}\!\left(\frac{U_s}{\Omega}\right)\cdot
\tau^r_j=\sum^n_{r=1}\frac{\nabla_{\!r}U_s}{\Omega}\,\tau^r_j
-\sum^n_{r=1}\sum^n_{m=1}\frac{p_m\,\nabla_{\!r}W^m\,U_s}{\Omega^2}
\,\tau^r_j.
\endgathered
\tag9.26
\endalign
$$
Right hand sides of the equalities \thetag{9.25} and \thetag{9.26}
are reflected in the above formula for $\theta_{ij}-\theta_{ji}$.
In deriving this formula for $\theta_{ij}-\theta_{ji}$ we have made
the following transformations for second term in right hand side of
commutator identity \thetag{9.19}:
$$
\gathered
\sum^n_{r=1}\sum^n_{\alpha=1}
\sum^n_{\gamma=1}\frac{p_\alpha}{\nu}\,D^{\alpha\gamma}_{sr}\,
\tau^r_j\,\xi_{\gamma i}=
\sum^n_{k=1}\sum^n_{m=1}
\sum^n_{r=1}\frac{p_m}{\nu}\,D^{mr}_{sk}\,
\tau^k_j\,\xi_{ri}=\\
=\sum^n_{k=1}\sum^n_{q=1}\sum^n_{m=1}\sum^n_{r=1}
\left(\,p_m\,D^{mr}_{sk}\,(P^q_r\,\nabla_{\!\boldsymbol\tau_i}n_q)
\,\tau^k_j-\frac{p_m\,p_q}{\nu\,\Omega}\,D^{mq}_{sk}\,U_r\,\tau^r_i\,
\tau^k_j\right).
\endgathered\quad
\tag9.27
$$
Terms from right hand side of \thetag{9.27} are also reflected
in the above formula for \linebreak $\theta_{ij}-\theta_{ji}$. Looking at
this formula, we see that compatibility equation \thetag{9.13}
providing compatibility of Pfaff equations \thetag{9.8} has the
following structure:
$$
\hskip -2em
\gathered
(\nu)^3\sum^n_{k=1}\sum^n_{q=1}\sum^n_{r=1}\sum^n_{s=1}
\frac{A^{rs}-A^{sr}}{\Omega}\,(P^q_r\,\nabla_{\!\boldsymbol
\tau_i}n_q)\,(P^k_s\,\nabla_{\!\boldsymbol\tau_j}n_k)\,+\\
+\,(\nu)^2\sum^n_{q=1}\sum^n_{r=1}\sum^n_{s=1}\frac{B^r_s}
{\Omega}\,(P^q_r\,\nabla_{\!\boldsymbol\tau_i}n_q)\,\tau^s_j\,
-(\nu)^2\sum^n_{q=1}\sum^n_{r=1}\sum^n_{s=1}\frac{B^r_s}{\Omega}
\,\times\\
\times\,(P^q_r\,\nabla_{\!\boldsymbol\tau_j}n_q)\,\tau^s_i
+\nu\,\sum^n_{r=1}\sum^n_{s=1}\frac{C_{rs}-C_{sr}}{\Omega}\,
\tau^r_i\,\tau^s_j=0.
\endgathered
\tag9.28
$$
Here $A^{rs}$, $B^r_s$, and $C_{rs}$ are components of three
extended tensor fields $\bold A$, $\bold B$, and $\bold C$
respectively. They are given by explicit formulas
$$
\align
&\hskip -2em
A^{rs}=\tilde\nabla^rW^s,
\tag9.29\\
&\hskip -2em
\aligned
B^r_s&=\tilde\nabla^rU_s+\sum^n_{m=1}\sum^n_{k=1}W^k\,p_m\,D^{mr}_{ks}\,-\\
&-\,\nabla_{\!s}W^r+\sum^n_{m=1}\frac{\tilde\nabla^mW^r-\tilde\nabla^rW^m}
{\Omega}\,U_s\,p_m,
\endaligned
\tag9.30\\
&\hskip -2em
\aligned
C_{rs}&=\nabla_{\!r}U_s
-\sum^n_{m=1}\frac{U_r\,\tilde\nabla^mU_s+U_s\,\nabla_{\!r}W^m}
{\Omega}\,p_m\,-
\\
&-\sum^n_{k=1}\sum^n_{q=1}\left(\,\shave{\sum^n_{m=1}}\frac{D^{mq}_{ks}
\,U_r}{\Omega}\,p_m+\frac{R^q_{krs}}{2}\right)W^k\,p_q.
\endaligned
\tag9.31
\endalign
$$
Further study of compatibility equations \thetag{9.28} with coefficients
\thetag{9.29}, \thetag{9.30}, \thetag{9.31} require some information
concerning geometry of hypersurfaces in non-metric geometry of manifolds
equipped with symmetric extended connection $\Gamma$ and generalized
Legendre transformation \thetag{2.2} given by extended vector field
$\bold V$ in $\bold p$-re\-presentation (see \thetag{2.4}).
\head
10. Geometry of hypersurfaces.
\endhead
    Let's fix some arbitrary point $q_0=(p_0,\bold p)$ of cotangent bundle
$T^*\!M$. This means that we fix some point $p_0\in M$ and some covector
$\bold p\in T^*_{p_0}(M)$. Assume that $\bold p\neq 0$. Then null-space
of covector $\bold p$ is a hyperplane in tangent space $T_{p_0}(M)$
(see \thetag{8.29}). Let $\sigma$ be some smooth hypersurface passing
through the point $p_0$ and tangent to null-space of fixed momentum
covector $\bold p$ at that point. If $\bold n=\bold n(p)$ is smooth
normal covector field on $\sigma$, then for $p=p_0$ we have
$$
\hskip -2em
\bold p=\nu_0\cdot\bold n(p_0)\text{, \ where \ }\nu_0\neq 0.
\tag10.1
$$
Taking constant $\nu_0\neq 0$ from \thetag{10.1}, we can expand it
up to a smooth nonzero function $\nu=\nu(p)$ on hypersurface $\sigma$
(or at least in some neighborhood of marked point $p_0$ on $\sigma$).
Function $\nu(p)$ satisfies the equality
$$
\hskip -2em
\nu(p_0)=\nu_0,
\tag10.2
$$
which is just the same as normalizing condition in \thetag{9.9}.
So, we can substitute $\nu(p)$ into \thetag{1.14} and use it for
defining shift of $\sigma$ along trajectories of Newtonian dynamical
system \thetag{7.2}. Now we assume that dynamical system \thetag{7.2}
satisfies additional normality condition (see definition~9.2). This
means that for any hypersurface $\sigma$ passing through our marked
point $p_0\in M$ and for any nonzero constant $\nu_0$ in
\thetag{10.2} Pfaff equations \thetag{9.8} are compatible. Hence
compatibility equations \thetag{9.28} are fulfilled. Note that
in \thetag{9.28} we have explicit entries of $\nu=\nu(p)$ and
implicit entries of $\nu$ through $\bold p=\nu\cdot\bold n$ in
arguments of extended tensor fields $\bold A$, $\bold B$, and
$\bold C$. Moreover, we have implicit entries of $\nu$ in covariant
derivatives
$$
\hskip -2em
\nabla_{\!\boldsymbol\tau_i}n_s=\frac{\partial n_s}{\partial y^i}
-\sum^n_{k=1}\sum^n_{r=1}\Gamma^k_{sr}\,n_k\,\tau^r_i
\tag10.3
$$
due to connection components $\Gamma^k_{sr}$ that depend on
momentum covector $\bold p=\nu\cdot\bold n$ (see \thetag{4.14}).
But in \thetag{9.28} we have no derivatives of function $\nu=
\nu(p)$. Therefore, if we write \thetag{9.28} only at our fixed
point $p=p_0$, we can replace all entries of $\nu$ by normalizing
constant $\nu_0$ from \thetag{10.2}. Using \thetag{10.1}, we can
express $\bold n=\bold n(p_0)$ through our fixed momentum covector
$\bold p$ at the point $p=p_0$:
$$
\hskip -2em
\bold n=\bold n(p_0)=\frac{\bold p}{\nu_0}.
\tag10.4
$$
In this form $\bold n(p_0)$ is not too specific property of hypersurface
$\sigma$. It determines only tangent hyperplane to $\sigma$ at fixed
point $p=p_0$. The only parameters in \thetag{9.28} that depend on fine
structure of hypersurface $\sigma$ at the point $p_0$ are covariant
derivatives \thetag{10.3}. Using them, in \cite{20} we have defined
a map $f\!: T_p(\sigma)\to T^*_p(M)$. Indeed, if $\bold\tau$ is some
arbitrary vector tangent to $\sigma$, then $\boldsymbol\tau=\alpha^1\cdot
\boldsymbol\tau_1+\ldots+\alpha^{n-1}\cdot\boldsymbol\tau_{n-1}$. Let
$$
\hskip -2em
f(\boldsymbol\tau)=\nabla_{\!\boldsymbol\tau}\bold n=
\sum^n_{r=1}\sum^{n-1}_{j=1}\left(\alpha^j\,
\nabla_{\!\boldsymbol\tau_j}n_r\right)
\kern -1pt\cdot dx^r.
\tag10.5
$$
It is easy to see that \thetag{10.5} defines linear map from tangent
hyperplane $T_p(\sigma)$ at the point $p\in\sigma$ to cotangent space
$T^*_p(M)$. We consider composite map
$$
\hskip -2em
\bold b=-\bold P^*\compos f\compos\bold P.
\tag10.6
$$
Projection operator $\bold P^*$ in \thetag{10.6} is a conjugate
operator for projector $\bold P$ with components \thetag{6.3}.
Remember that $\bold P$ projects onto the subspace $T_p(\sigma)
\in T_p(M)$. Therefore linear map $\bold b\!:T_p(M)\to T^*_p(M)$
is \pagebreak correctly defined by formula \thetag{10.6}.\par
    Linear map $\bold b$ defined by formula \thetag{10.6} is
associated with second fundamental form of hypersurface $\sigma$.
Indeed, let's define bilinear form
$$
\hskip -2em
b(\bold X,\bold Y)=\left<\bold b(\bold Y)\,|\,\bold X\right>.
\tag10.7
$$
Due to the presence of projection operators $\bold P$ and $\bold P^*$
in \thetag{10.6} we have
$$
\hskip -2em
b(\bold X,\bold Y)=b(\bold P(\bold X),\bold Y)=b(\bold X,
\bold P(\bold Y)).
\tag10.8
$$
\proclaim{Theorem 10.1} Bilinear form \thetag{10.7} defined by
linear map \thetag{10.6} is symmetric.
\endproclaim
When restricted to tangent space $T_p(\sigma)$ of hypersurface $\sigma$
bilinear form \thetag{10.7} yields second fundamental form of $\sigma$.
Its components
$$
\hskip -2em
\beta_{ij}=b(\boldsymbol\tau_i,\boldsymbol\tau_j)
\tag10.9
$$
define tensor field in inner geometry of hypersurface $\sigma$. From
\thetag{10.8} we derive the equality $b(\bold X,\bold Y)=b(\bold P
(\bold X),\bold P(\bold Y))$. It means that bilinear form
\thetag{10.7} and linear map \thetag{10.6} in outer space are
completely determined by components of second fundamental form
\thetag{10.9}. Further we need the following theorem.
\proclaim{Theorem 10.2} Let $q_0=(p_0,\bold p)$ be some fixed point
of cotangent bundle $T^*\!M$ with $\bold p\neq 0$ and let projector
$\bold P$ be the value of projector-valued extended tensor field
\thetag{6.3} at this point. Then any symmetric quadratic form $b$ in
$T_{p_0}(M)$ satisfying the equality \thetag{10.8} can be determined by
some hypersurface $\sigma$ passing through the point $p_0$ and tangent to
null-space of covector $\bold p$ at this point.
\endproclaim
   Theorems 10.1 and 10.2 are proved in paper \cite{20}. Though
these theorems are very important for further study of compatibility
equations \thetag{9.28}, we shall not repeat their proofs in present
paper.
\head
11. Additional normality equations.
\endhead
    As in previous section, let's fix some point $p_0\in M$ and some
covector $\bold p\neq 0$ at this point. This means that we fix some
point $q_0=(p_0,\bold p)$ of cotangent bundle $T^*\!M$. Let's fix
some arbitrary nonzero constant $\nu_0\neq 0$ and then use formula
\thetag{10.4} for to define another nonzero covector $\bold n\neq 0$
at our fixed point $p_0$. Further, let's consider various hypersurfaces
passing through the point $p_0$ tangent to null-space of covector
$\bold n$. For each such hypersurface $\sigma$ covector $\bold n$ is
normal covector at the point $p_0$. It can be expanded up to a smooth
normal covector field $\bold n=\bold n(p)$ (at least in some
neighborhood of marked point $p_0$). Therefore we can build $\sigma$
into a framework of shift construction defined by of Newtonian dynamical
\thetag{7.2}. If this dynamical system satisfies additional normality
condition (see definition~9.2), then we can choose smooth function
$\nu=\nu(p)$ normalized by the condition \thetag{10.2} and such that
compatibility equations \thetag{9.28} are fulfilled.\par
    Let $\boldsymbol\tau_1,\,\ldots,\,\boldsymbol\tau_{n-1}$ be basic
tangent vectors of $\sigma$. Applying linear map \thetag{10.6} to them,
we get a set of $n-1$ covectors $\boldsymbol\theta_1,\,\ldots,\,
\boldsymbol\theta_{n-1}$. In other words, we denote
$$
\pagebreak
\hskip -2em
\boldsymbol\theta_i=\bold b(\boldsymbol\tau_i).
\tag11.1
$$
Let's calculate components of covectors \thetag{11.1} in local chart.
Using formulas \thetag{10.5} and \thetag{10.6} defining linear map
$\bold b\!:T_p(M)\to T^*_p(M)$, we derive
$$
\hskip -2em
\theta_{ri}=-\sum^n_{q=1}P^q_r\,\nabla_{\!\boldsymbol\tau_i}n_q=
\sum^n_{s=1}b_{rs}\,\tau^s_i.
\tag11.2
$$
Here $b_{rs}$ are components of bilinear form \thetag{10.7} in local
chart. Comparing \thetag{9.28} and \thetag{11.2}, we see that
compatibility equations \thetag{9.28} can be written as follows:
$$
\hskip -2em
\gathered
(\nu_0)^3\sum^n_{k=1}\sum^n_{q=1}\sum^n_{r=1}\sum^n_{s=1}
\frac{A^{rs}-A^{sr}}{\Omega}\,(b_{rq}\,\tau^q_i)\,
(b_{sk}\,\tau^k_j)\,-\\
-\,(\nu_0)^2\sum^n_{q=1}\sum^n_{r=1}\sum^n_{s=1}\frac{B^r_s}
{\Omega}\,(b_{rq}\,\tau^q_i)\,\tau^s_j\,
+(\nu_0)^2\sum^n_{q=1}\sum^n_{r=1}\sum^n_{s=1}\frac{B^r_s}{\Omega}
\,\times\\
\times\,(b_{rq}\,\tau^q_j)\,\tau^s_i
+\nu_0\,\sum^n_{r=1}\sum^n_{s=1}\frac{C_{rs}-C_{sr}}{\Omega}\,
\tau^r_i\,\tau^s_j=0.
\endgathered
\tag11.3
$$
Varying hypersurface $\sigma$, we can vary components of bilinear
form $b$ in \thetag{11.3}. In particular, we can 1) change the sign
of $b$; 2) get zero quadratic form for $b$. Due to these facts
(they follow from theorem~10.2) compatibility equations \thetag{11.3}
split into three separate parts. Now they are written as follows:
$$
\align
&\hskip -2em
\sum^n_{r=1}\sum^n_{s=1}(A^{rs}-A^{sr})\,\theta_{ri}\,\theta_{sj}=0,
\tag11.4\\
&\hskip -2em
\sum^n_{r=1}\sum^n_{s=1}\sum^n_{q=1}(B^q_s\,b_{qr}-B^q_r\,b_{qs})\,
\tau^r_i\,\tau^s_j=0,
\tag11.5\\
&\hskip -2em
\sum^n_{r=1}\sum^n_{s=1}(C_{rs}-C_{sr})\,\tau^r_i\,\tau^s_j=0.
\tag11.6
\endalign
$$
Here in \thetag{11.4} we used \thetag{11.2} again. Vectors $\boldsymbol
\tau_1,\,\ldots,\,\boldsymbol\tau_{n-1}$ depend on the choice of local
chart on hypersurface $\sigma$. For a fixed point $p=p_0$ they can be
treated as arbitrary $n-1$ vectors forming base in tangent hyperplane
$T_p(\sigma)$. Similarly, due to theorem~10.2 covectors $\boldsymbol
\theta_1,\,\ldots,\,\boldsymbol\theta_{n-1}$ can be treated as arbitrary
$n-1$ covectors in null-space of vector $\bold W$ (see \thetag{8.32}).
Therefore \thetag{11.4} and \thetag{11.6} reduce to
$$
\align
&\hskip -2em
\sum^n_{r=1}\sum^n_{s=1}(A^{rs}-A^{sr})\,P^i_r\,P^j_s=0,
\tag11.7\\
&\hskip -2em
\sum^n_{r=1}\sum^n_{s=1}(C_{rs}-C_{sr})\,P^r_i\,P^s_j=0.
\tag11.8
\endalign
$$
In this form equations \thetag{11.7} and \thetag{11.8} do not depend on
any particular hypersurface we used to derive them.\par
    As for \thetag{11.5}, we should bring it to similar form independent
on $\sigma$. For this purpose remember that we can treat $\boldsymbol
\tau_1,\,\ldots,\,\boldsymbol\tau_{n-1}$ \pagebreak as arbitrary $n-1$
vectors in tangent hyperplane $T_p(\sigma)$. Therefore equation
\thetag{11.5} reduces to the following form:
$$
\hskip -2em
\sum^n_{r=1}\sum^n_{s=1}\sum^n_{q=1}(B^q_s\,b_{qr}-B^q_r\,b_{qs})
\,P^r_i\,P^s_j=0.
\tag11.9
$$
In the next step we rewrite \thetag{11.9} in coordinate-free form.
It looks like
$$
\hskip -2em
b(\bold P\compos\bold B\compos\bold P(\bold X),\bold Y)=
b(\bold X,\bold P\compos\bold B\compos\bold P(\bold Y))
\tag11.10
$$
Here $\bold B$ is linear operator in $T_{p_0}(M)$ determined by
components of extended tensor field \thetag{9.30}. Due to theorem~10.2
the above equality \thetag{11.10} means that composite operator
$\bold P\compos\bold B\compos\bold P$ is symmetric with respect to any
symmetric bilinear form $b$ in $T_{p_0}(M)$ for which \thetag{10.8} is
fulfilled. Hence composite operator $\bold P\compos\bold B\compos
\bold P$ can differ from projector $\bold P$ only by some scalar
factor $\lambda$:
$$
\hskip -2em
\bold P\compos\bold B\compos\bold P=\lambda\cdot\bold P.
\tag11.11
$$
This fact is proved in paper \cite{20}. Scalar factor $\lambda$ is
given by trace formula
$$
\hskip -2em
\lambda=\frac{\tr(\bold P\compos\bold B\compos\bold P)}{n-1}
=\frac{\tr(\bold B\compos\bold P)}{n-1}=
\sum^n_{r=1}\sum^n_{s=1}\frac{B^r_s\,P^s_r}{n-1}.
\tag11.12
$$
Formulas \thetag{11.11} and \thetag{11.12} written in local chart
yield required equations
$$
\hskip -2em
\sum^n_{r=1}\sum^n_{s=1}P^i_r\,B^r_s\,P^s_j=\sum^n_{r=1}\sum^n_{s=1}
\frac{B^r_s\,P^s_r}{n-1}\,P^i_j.
\tag11.13
$$
In form \thetag{11.13} equations \thetag{11.5} do not depend on any
particular hypersurface used in deriving these equations.\par
    Equations \thetag{11.7}, \thetag{11.8}, and \thetag{11.13} taken
together form a system of {\it additional normality equations}. Due
to the above notations \thetag{6.1}, \thetag{6.2}, \thetag{6.3},
\thetag{9.14}, \thetag{9.29}, \thetag{9.30}, and \thetag{9.31} they
are partial differential equations for components of extended vector
field $\bold V$ and extended covector field $\bold Q$ that determine
Newtonian dynamical system \thetag{7.2}. We have derived them assuming
that Newtonian dynamical system \thetag{7.2} satisfies additional
normality condition formulated in definition~9.2. Conversely, if
additional normality equations are fulfilled, then compatibility
equations \thetag{9.28} turn to identities. Therefore Pfaff equations
\thetag{9.8} appear to be compatible. Thus, we have proved the
following theorem analogous to theorem~8.1.
\proclaim{Theorem 11.1} Additional normality condition for Newtonian
dynamical system \thetag{7.2} is equivalent to the system of additional
normality equations \thetag{11.7}, \thetag{11.8}, and \thetag{11.13}
that should be fulfilled at all points $q=(p,\bold p)$ of cotangent
bundle $T^*\!M$, where $\bold p\neq 0$.
\endproclaim
    Now suppose that both weak and additional normality conditions
are fulfilled. In this case all deviation functions \thetag{9.3} satisfy
second order ordinary differential equation \thetag{8.8} and we can
provide initial data \thetag{9.5} for them by proper choice of function
$\nu=\nu(p)$ on any predefined hypersurface $\sigma$ in $M$. \pagebreak 
Then all deviation functions do vanish, and we have normal shift of
hypersurface $\sigma$. So, we see that weak and additional normality
conditions are complementary to each other, and if both are fulfilled,
we can arrange normal shift of any predefined hypersurface along
trajectories of Newtonian dynamical system \thetag{7.2}.
\definition{Definition 11.1} We say that Newtonian dynamical system
satisfies {\bf complete} normality condition if both {\bf weak} and
{\bf additional} normality conditions for this system are fulfilled.
\enddefinition
   According to theorems~8.1 and 11.1, weak and additional normality
conditions are equivalent to weak and additional normality equations
for parameters $\bold V$ and $\bold Q$ determining Newtonian dynamical
system \thetag{7.2}. Therefore complete normality condition is
equivalent to complete system of normality equations including
\thetag{8.21}, \thetag{11.7}, \thetag{11.8}, and \thetag{11.13}. As
we noted just above, complete normality condition is sufficient for
strong normality condition to be fulfilled (see definition~1.5). We
shall strengthen this result in the next section.
\head
12. Equivalence of strong and complete normality conditions.
\endhead
    Part of the statement declared in the title of this section is
already proved. Indeed, we know that complete normality condition
implies strong normality condition. Let's prove converse implication.
Assuming that Newtonian dynamical system \thetag{7.2} satisfies
strong normality condition, we should prove that it satisfies weak
and additional normality conditions.\par
     In the first step let's prove that additional normality condition
is fulfilled. For this purpose let's take some arbitrary hypersurface
$\sigma$ with marked point $p=p_0$ and smooth normal covector field
$\bold n=\bold n(p)$ in some neighborhood of marked point. Then let's
take some nonzero constant $\nu_0\neq 0$ and let's apply strong normality
condition, which is fulfilled by assumption (see definition~1.5). As a
result we get smooth function $\nu=\nu(p)$ normalized by condition
\thetag{10.2} and such that it provide initial data \thetag{1.14} for
normal shift of $\sigma$. Due to normality of shift all deviation functions
\thetag{9.3} are identically zero. Hence initial conditions \thetag{9.5}
for them are fulfilled. Writing \thetag{9.5} explicitly, we find that
our function $\nu=\nu(p)$ is a solution for Pfaff equations \thetag{9.8}.
Now, varying constant $\nu_0\neq 0$ in \thetag{10.2}, we prove that 
Pfaff equations \thetag{9.8} are compatible (see definition~9.1). Thus,
additional normality condition is proved (see definition~9.2).\par
    In the second step we shall derive weak normality condition assuming
that strong normality condition is fulfilled. This is a little bit more
complicated. For this purpose we fix some point $q_0=(p_0,\bold p)$ of
cotangent bundle $T^*\!M$ with $\bold p\neq 0$. Initial point $p_0$ and
momentum covector $\bold p$ at this point form initial data for Newtonian
dynamical system \thetag{7.2}. They define a trajectory $p=p(t)$ passing
through initial point $p_0$. Null-space of initial covector $\bold p$ is
a hyperplane in tangent space $T_{p_0}(M)$. Let's consider various
hypersurfaces passing through initial point $p_0$ tangent to this
hyperplane and denote by $\sigma$ one of them. If $\bold n=\bold n(p)$ is
normal covector of this hypersurface $\sigma$, then at the point $p=p_0$
we have the equality \thetag{10.1} that determine normalizing constant
$\nu_0\neq 0$ for \thetag{10.2}. Applying strong normality condition (see
definition~1.5), we find smooth function $\nu=\nu(p)$ on $\sigma$ in some
neighborhood of initial point $p_0$ that provide initial data \thetag{1.14}
for normal shift of $\sigma$. Thus, our fixed trajectory $p=p(t)$ passing
through initial point $p_0$ appears to be shift trajectory among many
others. Due to normality of shift all corresponding deviation functions
\thetag{9.3} on this trajectory are identically zero. Therefore we can write
$$
\hskip -2em
\ddot\varphi_i\,\hbox{\vrule height 8pt depth 8pt width 0.5pt}_{\,t=0}=0.
\tag12.1
$$
Note that initial conditions \thetag{9.5} are also fulfilled. Therefore
we can use formula \thetag{8.20} for $\ddot\varphi$ in left hand side
of \thetag{12.1}. Here it is written as follows:
$$
\hskip -2em
\ddot\varphi_i\,\hbox{\vrule height 8pt depth 8pt width 0.5pt}_{\,t=0}
=\sum^n_{k=1}\sum^n_{r=1}\alpha^r\,P^k_r\cdot\xi_{ki}\,\hbox{\vrule
height 8pt depth 8pt width 0.5pt}_{\,t=0}+\sum^n_{k=1}\sum^n_{r=1}\eta_r
\,P^r_k\cdot\tau^k_i\,\hbox{\vrule height 8pt depth 8pt width
0.5pt}_{\,t=0}.
\tag12.2
$$
Vectors $\boldsymbol\tau_1,\,\ldots,\,\boldsymbol\tau_{n-1}$ at
initial instant of time $t=0$ form base in tangent hyperplane to
initial hypersurface $\sigma$. As we noted in section~11, for
fixed point $p=p_0$ they can be treated as arbitrary $n-1$ vectors
in null-space of momentum covector $\bold p=\nu_0\cdot\bold n(p_0)$.
For components of covectors $\boldsymbol\xi_1,\,\ldots,\,\boldsymbol
\xi_{n-1}$ at initial instant of time $t=0$ we can use formula
\thetag{9.20}. Applying \thetag{9.20} and \thetag{11.2}, we get
$$
\hskip -2em
\sum^n_{k=1}P^k_r\cdot\xi_{ki}\,\hbox{\vrule
height 8pt depth 8pt width 0.5pt}_{\,t=0}=
\sum^n_{s=1}\nu_0\,P^s_r\,\nabla_{\!\boldsymbol\tau_i}n_s=
-\sum^n_{s=1}\nu_0\,b_{rs}\,\tau^s_i.
\tag12.3
$$
Combining \thetag{12.1}, \thetag{12.2}, and \thetag{12.3}, we derive
the equality
$$
\hskip -2em
-\sum^n_{k=1}\sum^n_{r=1}\sum^n_{s=1}\nu_0\,\alpha^r\,P^k_r
b_{ks}\,\tau^s_i+\sum^n_{k=1}\sum^n_{r=1}\eta_r
\,P^r_k\cdot\tau^k_i=0.
\tag12.4
$$
Just like in \thetag{11.3}, by varying hypersurface $\sigma$ and by
applying theorem~10.2 we can break \thetag{12.4} into two separate
equalities
$$
\xalignat 2
&\hskip -2em
\sum^n_{k=1}\sum^n_{r=1}\alpha^r\,P^k_r\cdot\theta_{ki}=0,
&&\sum^n_{k=1}\sum^n_{r=1}\eta_r\,P^r_k\cdot\tau^k_i=0,
\endxalignat
$$
which are equivalent to weak normality equations \thetag{8.21}.
Applying theorem~8.1, we find that weak normality condition is
fulfilled, i\.\,e\. strong normality condition implies weak
normality condition. Ultimately, we have proved the following
theorem.
\proclaim{Theorem 12.1} Strong and complete normality conditions
for Newtonian dynamical system \thetag{7.2} are equivalent to each
other.
\endproclaim
\head
13. Connection invariance.
\endhead
    Theorem~12.1 is a basic result in the theory of dynamical Newtonian
dynamical systems admitting normal shift, while definition~1.5 is basic
definition of this theory. Comparing them we see that strong normality
condition formulated in definition~1.5 is applicable either to general
Newtonian dynamical system of the form \thetag{2.3}, and to special one
given by the equations \thetag{7.2}. Theorem~12.1 is formulated only for
Newtonian dynamical system \thetag{7.2}, which implies presence of some
symmetric extended connection $\Gamma$ in $M$. Theorem~5.1 gives one
way to avoid this discrepancy. We can use symmetric affine connection
\thetag{5.8} canonically associated with dynamical system \thetag{2.3}
and by means of this connection we can rewrite \thetag{2.3} in form of
\thetag{7.2} (see formula \thetag{7.1}). However, there is another way.
Below we shall prove that the whole theory constructed in sections~6--13
is invariant under gauge transformations changing one connection for
another:
$$
\hskip -2em
\Gamma^k_{ij}\to \Gamma^k_{ij}+T^k_{ij}.
\tag13.1
$$
Here $T^k_{ij}$ are components of some symmetric extended tensor field
$\bold T$ of type $(1,2)$. Applying gauge transformation \thetag{13.1}
to dynamical system \thetag{7.2}, we change covariant derivatives
$\nabla_{\!t}p_i$ in left hand side. In order to keep corresponding
connection-free equations \thetag{2.3} unchanged we should change
components of covector $\bold Q$ as follows:
$$
\hskip -2em
Q_i\to Q_i-\sum^n_{k=1}\sum^n_{s=1}T^k_{is}\,p_k\,V^s.
\tag13.2
$$\par
    Having fixed gauge transformations by formulas \thetag{13.1} and
\thetag{13.2}, now we shall apply them to all normality equations
\thetag{8.21}, \thetag{11.7}, \thetag{11.8}, and \thetag{11.13} for
to prove their invariance under these transformations. First of all
note that vector field $\bold V$ in \thetag{7.2}, vector field
$\bold W$ introduced by formula \thetag{6.1}, scalar field $\Omega$
given by formula \thetag{6.2}, and projector field $\bold P$ with
components \thetag{6.3} are invariant under gauge transformations
defined by \thetag{13.1}, \thetag{13.2}:
$$
\xalignat 2
&\hskip -2em
V^s\to V^s,&&W^s\to W^s,
\tag13.3\\
&\hskip -2em
\Omega\to \Omega,&&P^i_j\to P^i_j.
\tag13.4
\endxalignat
$$
As for covector field $\bold U$ in \thetag{9.14}, here we have
the following transformation rule:
$$
\hskip -2em
U_s\to U_s+\sum^n_{q=1}\sum^n_{m=1}p_m\,T^m_{sq}\,W^q.
\tag13.5
$$
Applying \thetag{13.1} to curvature tensors \thetag{8.6} and
\thetag{8.7}, we derive
$$
\align
&\hskip -4em
\aligned
&R^k_{rij}\to R^k_{rij}+\nabla_{\!i}T^k_{jr}-\nabla_{\!j}T^k_{ir}
-\sum^n_{s=1}\sum^n_{m=1}p_s\,D^{km}_{jr}\,T^s_{mi}\,+\\
&+\,\sum^n_{s=1}\sum^n_{m=1}p_s\,D^{km}_{ir}\,T^s_{mj}
+\sum^n_{m=1}\left(T^k_{im}\,T^m_{jr}-T^k_{jm}\,T^m_{ir}\right)+\\
&+\,\sum^n_{s=1}\sum^n_{m=1}p_s\,T^s_{mi}\,\tilde\nabla^mT^k_{jr}
-\sum^n_{s=1}\sum^n_{m=1}p_s\,T^s_{mj}\,\tilde\nabla^mT^k_{ir},
\endaligned\hskip -1em
\tag13.6\\
\vspace{2ex}
&\hskip -4em
D^{kr}_{ij}\to D^{kr}_{ij}-\tilde\nabla^rT^k_{ij}.
\tag13.7
\endalign
$$
Now we can apply \thetag{13.1}, \thetag{13.2}, \thetag{13.3},
\thetag{13.4}, and \thetag{13.7} to vector field $\boldsymbol
\alpha$ with components \thetag{8.12} used in weak normality
equations \thetag{8.21}:
$$
\allowdisplaybreaks
\gather
\alpha^k\to\alpha^k-\sum^n_{i=1}\sum^n_{r=1}\sum^n_{s=1}\tilde\nabla^kV^i
\,T^r_{is}\,p_r\,V^s+\left(\,\shave{\sum^n_{r=1}\sum^n_{s=1}\sum^n_{i=1}}
T^k_{rs}\,\tilde\nabla^sV^i\,p_i\,V^r\,+\right.\\
\left.+\,\shave{\sum^n_{r=1}\sum^n_{s=1}\sum^n_{i=1}}T^i_{rs}\,\tilde
\nabla^kV^s\,p_i\,V^r+\shave{\sum^n_{r=1}\sum^n_{s=1}\sum^n_{i=1}
\sum^n_{m=1}}p_m\,T^m_{rs}\,\tilde\nabla^s\tilde\nabla^kV^i\,p_i\,
V^r\!\right)-\\
-\,\sum^n_{r=1}\sum^n_{s=1}\sum^n_{i=1}\sum^n_{m=1}\tilde\nabla^r\tilde
\nabla^kV^i\,p_i\,T^m_{rs}\,p_m\,V^s-\sum^n_{r=1}\sum^n_{i=1}
\left(\,\shave{\sum^n_{s=1}}\tilde\nabla^rV^i\,p_i\,T^k_{rs}\,V^s\,
-\right.\\
\left.-\,\sum^n_{m=1}\sum^n_{s=1}\tilde\nabla^rV^i\,p_i\,\tilde\nabla^k
T^m_{rs}\,p_m\,V^s-\sum^n_{m=1}\sum^n_{s=1}\tilde\nabla^rV^i\,p_i
\,T^m_{rs}\,p_m\,\tilde\nabla^kV^s\!\right)+\\
+\,\sum^n_{r=1}\sum^n_{s=1}\sum^n_{i=1}\sum^n_{m=1}\tilde\nabla^k
T^s_{rm}\,\tilde\nabla^rV^i\,p_i\,p_s\,V^m+\sum^n_{r=1}\sum^n_{i=1}
\left(\,\shave{\sum^n_{s=1}}T^i_{rs}\,V^s\,p_i\tilde\nabla^kV^r\,
+\right.\\
\left.+\,\shave{\sum^n_{s=1}\sum^n_{m=1}}p_m\,T^m_{rs}\,\tilde\nabla^s
V^i\,p_i\tilde\nabla^kV^r\!\right)-\sum^n_{r=1}\sum^n_{s=1}\sum^n_{m=1}
\tilde\nabla^kV^r\,T^m_{rs}\,p_m\,V^s.
\endgather
$$
Looking attentively at the above formula, we see that almost all
terms in right hand side do cancel each other. As a result we get
the following transformation rule:
$$
\hskip -2em
\alpha^k\to\alpha^k.
\tag13.8
$$
Formula \thetag{8.13} for components of extended covector field
$\boldsymbol\beta$ is more complicated than formula \thetag{8.12}.
Let's simplify it using notations \thetag{9.14}:
$$
\hskip -2em
\gathered
\beta_k=\sum^n_{r=1}\nabla_{\!r}U_k\,V^r+\sum^n_{r=1}
\tilde\nabla^rU_k\,Q_r+\sum^n_{r=1}\nabla_{\!k}V^r\,U_r\,+\\
+\sum^n_{r=1}\nabla_{\!k}Q_r\,W^r-\sum^n_{r=1}\sum^n_{s=1}
\sum^n_{m=1}\left(R^s_{rmk}\,V^m-D^{sm}_{rk}\,Q_m\right)W^r\,p_s.
\endgathered
\tag13.9
$$
Formula \thetag{13.9} for $\beta^k$ is still rather complicated. 
Therefore we perform some preliminary calculations. Using formula
\thetag{13.5}, we derive
$$
\align
&\gathered
\nabla_{\!r}U_s\to\nabla_{\!r}U_s
+\sum^n_{q=1}\sum^n_{m=1}p_m\,\nabla_{\!r}T^m_{sq}\,\,W^q
+\sum^n_{q=1}\sum^n_{m=1}p_m\,T^m_{sq}\,\nabla_{\!r}W^q\,-\\
-\,\sum^n_{k=1}T^k_{rs}\,U_k-\sum^n_{k=1}\sum^n_{q=1}\sum^n_{m=1}
p_m\,T^k_{rs}\,\,T^m_{kq}\,\,W^q\,+\sum^n_{k=1}\sum^n_{u=1}p_u\,
T^u_{rk}\,\times\\
\times\,\tilde\nabla^kU_s+\sum^n_{k=1}\sum^n_{u=1}
\sum^n_{q=1}p_u\,T^u_{rk}\,T^k_{sq}\,W^q+\sum^n_{k=1}\sum^n_{u=1}
\sum^n_{q=1}\sum^n_{m=1}p_u\,p_m\,\times\\
\times\,T^u_{rk}\,\tilde\nabla^kT^m_{sq}\,W^q+\sum^n_{k=1}
\sum^n_{u=1}\sum^n_{q=1}\sum^n_{m=1}p_u\,p_m\,T^u_{rk}\,T^m_{sq}
\,\tilde\nabla^kW^q,
\endgathered\qquad
\tag13.10\\
&\gathered
\tilde\nabla^rU_s\to\tilde\nabla^rU_s
+\sum^n_{q=1}T^r_{sq}\,W^q+\sum^n_{q=1}\sum^n_{m=1}p_m\,\tilde
\nabla^rT^m_{sq}\,W^q\,+\\
+\,\sum^n_{q=1}\sum^n_{m=1}p_m\,T^m_{sq}\,\tilde\nabla^rW^q.
\endgathered
\tag13.11
\endalign
$$
For covariant derivative $\nabla_{\!k}V^i$ we have the following
transformation rule:
$$
\hskip -2em
\nabla_{\!k}V^i\to\nabla_{\!k}V^i+\sum^n_{s=1}T^i_{ks}\,V^s
+\sum^n_{s=1}\sum^n_{m=1}p_m\,T^m_{ks}\,\tilde\nabla^sV^i.
\tag13.12
$$
Combining formulas \thetag{13.11} and \thetag{13.2}, we derive
$$
\gathered
\sum^n_{r=1}\tilde\nabla^rU_k\,Q_r\to\sum^n_{r=1}\tilde\nabla^rU_k
\,Q_r+\sum^n_{r=1}\sum^n_{q=1}T^r_{kq}\,W^q\,Q_r+\sum^n_{r=1}
\sum^n_{q=1}\sum^n_{m=1}p_m\,\times\\
\times\,\tilde\nabla^rT^m_{kq}\,W^q\,Q_r+\sum^n_{r=1}\sum^n_{q=1}
\sum^n_{m=1}p_m\,T^m_{kq}\,\tilde\nabla^rW^q\,Q_r-\sum^n_{r=1}
\sum^n_{s=1}\sum^n_{u=1}\tilde\nabla^rU_k\,\times\\
\times\,T^u_{rs}\,p_u\,V^s-\sum^n_{r=1}\sum^n_{s=1}\sum^n_{q=1}
\sum^n_{u=1}T^r_{kq}\,W^q\,T^u_{rs}\,p_u\,V^s
-\sum^n_{r=1}\sum^n_{s=1}\sum^n_{q=1}\sum^n_{u=1}\sum^n_{m=1}
p_u\,\times\\
\times\,p_m\,\tilde\nabla^rT^m_{kq}\,W^q\,T^u_{rs}\,V^s
-\sum^n_{r=1}\sum^n_{s=1}\sum^n_{q=1}\sum^n_{u=1}\sum^n_{m=1}
p_u\,p_m\,T^m_{kq}\,T^u_{rs}\,\tilde\nabla^rW^q\,V^s.
\endgathered
$$
In a similar way, combining formulas \thetag{13.12} and \thetag{13.5},
we derive
$$
\gather
\sum^n_{r=1}\nabla_{\!k}V^r\,U_r\to\sum^n_{r=1}\nabla_{\!k}V^r\,U_r
+\sum^n_{r=1}\sum^n_{s=1}T^r_{ks}\,V^s\,U_r+\sum^n_{r=1}\sum^n_{s=1}
\sum^n_{m=1}p_m\,\times\\
\times\,T^m_{ks}\,\tilde\nabla^sV^r\,U_r+\sum^n_{r=1}\sum^n_{q=1}
\sum^n_{m=1}p_m\,T^m_{rq}\,\nabla_{\!k}V^r\,W^q+\sum^n_{r=1}
\sum^n_{s=1}\sum^n_{q=1}\sum^n_{m=1}p_m\,\times\\
\times\,T^r_{ks}\,T^m_{rq}\,W^q\,V^s+\sum^n_{r=1}\sum^n_{s=1}
\sum^n_{q=1}\sum^n_{u=1}\sum^n_{m=1}p_u\,p_m\,T^u_{rq}\,T^m_{ks}
\,\tilde\nabla^sV^r\,W^q.
\endgather
$$
Transformation rule for $\nabla_{\!k}Q_r$ is similar to \thetag{13.10}:
$$
\gathered
\nabla_{\!k}Q_r\to\nabla_{\!k}Q_r-\sum^n_{s=1}\sum^n_{u=1}
p_u\,\nabla_{\!k}T^u_{rs}\,V^s-\sum^n_{s=1}\sum^n_{u=1}
p_u\,T^u_{rs}\,\nabla_{\!k}V^s\,-\\
-\sum^n_{q=1}T^q_{kr}\,Q_q
+\sum^n_{s=1}\sum^n_{u=1}\sum^n_{q=1}p_u\,T^q_{kr}\,T^u_{rs}\,V^s
+\sum^n_{m=1}\sum^n_{q=1}p_m\,T^m_{kq}\,\times\\
\times\,\tilde\nabla^qQ_r-\sum^n_{s=1}\sum^n_{q=1}\sum^n_{m=1}p_m\,
T^m_{kq}\,T^q_{rs}\,V^s-\sum^n_{s=1}\sum^n_{q=1}\sum^n_{m=1}
\sum^n_{u=1}p_u\,p_m\,\times\\
\times\,T^m_{kq}\,\tilde\nabla^qT^u_{rs}\,V^s-\sum^n_{s=1}\sum^n_{q=1}
\sum^n_{m=1}\sum^n_{u=1}p_u\,p_m\,T^m_{kq}\,T^u_{rs}\,\tilde\nabla^qV^s.
\endgathered\quad
\tag13.13
$$
For term with curvature tensor $D^{sm}_{rk}$ in \thetag{13.9} we use
transformation rule \thetag{13.7}:
$$
\allowdisplaybreaks
\gather
\sum^n_{r=1}\sum^n_{s=1}\sum^n_{m=1}p_m\,D^{ms}_{rk}\,Q_s\,W^r\to
\sum^n_{r=1}\sum^n_{s=1}\sum^n_{m=1}p_m\,D^{ms}_{rk}\,Q_s\,W^r\,-\\
-\,\sum^n_{r=1}\sum^n_{s=1}\sum^n_{m=1}p_m\,\tilde\nabla^sT^m_{rk}
\,Q_s\,W^r-\sum^n_{r=1}\sum^n_{s=1}\sum^n_{q=1}\sum^n_{u=1}
\sum^n_{m=1}p_m\,D^{ms}_{rk}\,T^u_{sq}\,p_u\,V^q\,W^r\,+\\
+\,\sum^n_{r=1}\sum^n_{s=1}\sum^n_{q=1}\sum^n_{u=1}\sum^n_{m=1}
p_m\,\tilde\nabla^sT^m_{rk}\,T^u_{sq}\,p_u\,V^q\,W^r.
\endgather
$$
And finally, for term with another curvature tensor $R^s_{rmk}$ in
formula \thetag{13.9} we use transformation rule \thetag{13.6}, which
is more complicated:
$$
\gathered
\sum^n_{r=1}\sum^n_{s=1}\sum^n_{m=1}p_m\,R^m_{rsk}\,V^s\,W^r
\to\sum^n_{r=1}\sum^n_{s=1}\sum^n_{m=1} p_m\,R^m_{rsk}\,V^s\,W^r\,+\\
+\,\sum^n_{r=1}\sum^n_{s=1}\sum^n_{m=1}p_m\,\nabla_{\!s}T^m_{kr}
\,V^s\,W^r-\,\sum^n_{r=1}\sum^n_{s=1}\sum^n_{m=1}p_m\,\nabla_{\!k}
T^m_{sr}\,V^s\,W^r\,-\\
-\,\sum^n_{r=1}\sum^n_{s=1}\sum^n_{m=1}\sum^n_{u=1}\sum^n_{q=1}
p_m\,p_u\left(D^{mq}_{kr}\,T^u_{qs}\,V^s\,W^r
-D^{mq}_{sr}\,T^u_{qk}\,V^s\,W^r\right)+\\
+\,\sum^n_{r=1}\sum^n_{s=1}\sum^n_{m=1}\sum^n_{q=1}\left(\,p_m\,T^m_{sq}
\,T^q_{kr}\,W^r\,V^s-p_m\,T^m_{kq}\,T^q_{sr}\,W^r\,V^s\right)+\\
+\,\sum^n_{r=1}\sum^n_{s=1}\sum^n_{m=1}\sum^n_{q=1}\sum^n_{u=1}
p_u\,p_m\left(T^u_{qs}\,\tilde\nabla^qT^m_{kr}-T^u_{qk}\,\tilde\nabla^q
T^m_{sr}\right)W^r\,V^s.
\endaligned
$$
Combining \thetag{13.10}, \thetag{13.11}, \thetag{13.12}, \thetag{13.13}
and other above formulas, for components of covector field $\boldsymbol
\beta$ we derive the following very simple transformation rule:
$$
\hskip -2em
\beta_k\to\beta_k+\sum^n_{m=1}\sum^n_{q=1}T^m_{kq}\,p_m\,\alpha^q.
\tag13.14
$$
In deriving \thetag{13.14} we used simplified version of formula
\thetag{8.12} for vector field $\boldsymbol\alpha$:
$$
\gathered
\alpha^k=\sum^n_{r=1}\tilde\nabla^kV^r\,U_r+\sum^n_{r=1}\nabla_{\!r}
W^k\,V^r+\sum^n_{r=1}\tilde\nabla^rW^k\,Q_r\,+\\
+\,\sum^n_{r=1}W^r\,\tilde\nabla^kQ_r-\sum^n_{r=1}\sum^n_{s=1}
\sum^n_{q=1}p_s\,D^{sk}_{rq}\,W^r\,V^q.
\endgathered
$$
In a similar way, using notations \thetag{9.14}, we can simplify formula
\thetag{8.18}:
$$
\hskip -2em
\eta_k=\beta_k-\sum^n_{s=1}\frac{U_k\,\alpha^s\,p_s}{\Omega}.
\tag13.15
$$
Substituting \thetag{13.14} into \thetag{13.15} and using \thetag{13.4},
\thetag{13.5}, and \thetag{13.8}, we obtain
$$
\pagebreak
\hskip -2em
\eta_k\to\eta_k+\sum^n_{s=1}\sum^n_{m=1}\sum^n_{q=1}
p_m\,T^m_{kq}\,P^q_s\,\alpha^s.
\tag13.16
$$
Formula \thetag{13.16} means that in general extended covector
field $\boldsymbol\eta$ is not invariant under gauge transformations
\thetag{13.1}, \thetag{13.2}. However, if first equation in
\thetag{8.21} is fulfilled, then formula \thetag{13.16} simplifies.
It turns to
$$
\eta_k\to\eta_k.
$$
\proclaim{Theorem 13.1} Weak normality equations \thetag{8.21} for
Newtonian dynamical system \thetag{7.2} are invariant under gauge
transformations \thetag{13.1}, \thetag{13.2}.
\endproclaim
    Now let's proceed with our calculations for additional normality
equations \thetag{11.7}, \thetag{11.8}, \thetag{11.13}. Scalar field
$\bold A$ given by formula \thetag{9.29} does not depend on $\bold Q$
and on connection components $\Gamma^k_{ij}$. Therefore we have 
$$
\hskip -2em
A_{rs}\to A_{rs}.
\tag13.17
$$
For $B^r_s$, using \thetag{13.5}, \thetag{13.7}, and \thetag{13.11},
from \thetag{9.30} we derive
$$
\hskip -2em
B^r_s\to B^r_s-\sum^n_{m=1}\sum^n_{q=1}\sum^n_{k=1}
p_k\,T^k_{sq}\left(\tilde\nabla^mW^r-\tilde\nabla^rW^m\right)P^q_m.
\tag13.18
$$
Now let's substitute \thetag{13.18} into \thetag{11.13} and use
\thetag{9.29}. As a result we get
$$
\hskip -2em
\gathered
\sum^n_{r=1}\sum^n_{s=1}P^i_r\,B^r_s\,P^s_j\to\sum^n_{r=1}
\sum^n_{s=1}P^i_r\,B^r_s\,P^s_j\,+\\
+\sum^n_{r=1}\sum^n_{s=1}\sum^n_{m=1}\sum^n_{q=1}\sum^n_{k=1}
p_k\,T^k_{sq}\,P^s_j\left(A^{rm}-A^{mr}\right)P^i_r\,P^q_m.
\endgathered
\tag13.19
$$
Looking at \thetag{13.19}, we see that in general left hand
side of the equation \thetag{11.13} is not invariant under
gauge transformations \thetag{13.1}, \thetag{13.2}. However,
if equations \thetag{11.7} are fulfilled, then formula
\thetag{13.19} reduces to
$$
\hskip -2em
\sum^n_{r=1}\sum^n_{s=1}P^i_r\,B^r_s\,P^s_j\to\sum^n_{r=1}
\sum^n_{s=1}P^i_r\,B^r_s\,P^s_j.
\tag13.20
$$
For the expression in right hand side of \thetag{11.13} under the
same assumption we have
$$
\hskip -2em
\sum^n_{r=1}\sum^n_{s=1}\frac{B^r_s\,P^s_r}{n-1}\,P^i_j\to
\sum^n_{r=1}\sum^n_{s=1}\frac{B^r_s\,P^s_r}{n-1}\,P^i_j.
\tag13.21
$$
This follows from trace formula \thetag{11.12} for scalar factor
$\lambda$ in \thetag{11.11}.\par
    Now let's apply gauge transformation \thetag{13.1}, \thetag{13.2}
to tensor field $\bold C$ with components \thetag{9.31}. By rather
huge, but direct calculations we find:
$$
\allowdisplaybreaks
\gather
C_{rs}\to C_{rs}+\ldots+\sum^n_{k=1}\sum^n_{q=1}\sum^n_{u=1}p_u\,T^u_{rk}
\,P^k_q\left(\tilde\nabla^qU_s-\nabla_sW^q\vphantom{\shave{\sum^n_{m=1}
\sum^n_{v=1}}}\,+\right.\\
\left.+\,\shave{\sum^n_{m=1}\sum^n_{v=1}}D^{vq}_{ms}\,W^m\,p_v\!\right)
+\sum^n_{k=1}\sum^n_{u=1}\sum^n_{m=1}p_u\,T^u_{rk}\frac{\tilde\nabla^m
W^k-\tilde\nabla^kW^m}{\Omega}\,U_s\,p_m\,+\\
+\,\sum^n_{k=1}\sum^n_{q=1}\sum^n_{u=1}\sum^n_{m=1}\sum^n_{v=1}
p_u\,T^u_{rk}\left(P^k_q\,\tilde\nabla^qW^m\,p_v\,T^v_{sm}
-\tilde\nabla^kW^m\,\frac{W^q\,p_m}{\Omega}\,p_v
\,T^v_{sq}\right).
\endgather
$$
By dots we denoted terms symmetric with respect to indices $r$ and $s$.
They do not affect ultimate equation \thetag{11.8}, therefore we need
not keep them in explicit form. In right hand side of the above formula
we see three distinct terms with sums. For the beginning let's transform
second term using formula \thetag{8.16}:
$$
\gather
\sum^n_{k=1}\sum^n_{u=1}\sum^n_{m=1}p_u\,T^u_{rk}\,\frac{\tilde\nabla^m
W^k-\tilde\nabla^kW^m}{\Omega}\,U_s\,p_m=\sum^n_{k=1}\sum^n_{u=1}
\sum^n_{m=1}\sum^n_{q=1}p_u\,T^u_{rk}\,P^k_q\,\times\\
\times\,\frac{\tilde\nabla^mW^q-\tilde\nabla^qW^m}{\Omega}\,U_s\,p_m
+\sum^n_{u=1}\sum^n_{k=1}\left(\,\shave{\sum^n_{q=1}\sum^n_{m=1}}
p_q\,\frac{\tilde\nabla^mW^q-\tilde\nabla^qW^m}{\Omega}\,p_m\!\right)
\times\\
p_u\,T^u_{rq}\,\frac{W^q}{\Omega}\,U_s=\sum^n_{k=1}\sum^n_{u=1}
\sum^n_{m=1}\sum^n_{q=1}p_u\,T^u_{rk}\,P^k_q\,\frac{\tilde\nabla^mW^q
-\tilde\nabla^qW^m}{\Omega}\,U_s\,p_m.
\endgather
$$
Substituting this result into formula for $C_{rs}$ and taking into
account \thetag{9.30}, we get
$$
\gather
C_{rs}\to C_{rs}+\ldots+\sum^n_{k=1}\sum^n_{q=1}\sum^n_{u=1}p_u\,T^u_{rk}
\,P^k_q\,B^q_s+\sum^n_{k=1}\sum^n_{q=1}\sum^n_{u=1}\sum^n_{m=1}\sum^n_{v=1}
p_u\,T^u_{rk}\,\times\\
\times\,P^k_q\,\tilde\nabla^qW^m\,p_v\,T^v_{sm}
-\sum^n_{k=1}\sum^n_{q=1}\sum^n_{u=1}\sum^n_{m=1}\sum^n_{v=1}
p_u\,T^u_{rk}\,\tilde\nabla^kW^m\,\frac{W^q\,p_m}{\Omega}\,p_v
\,T^v_{sq}.
\endgather
$$
Now we apply the same trick with formula \thetag{8.16} to last term
in the above formula:
$$
\gathered
\sum^n_{k=1}\sum^n_{q=1}\sum^n_{u=1}\sum^n_{m=1}\sum^n_{v=1}
p_u\,T^u_{rk}\,\tilde\nabla^kW^m\,\frac{W^q\,p_m}{\Omega}\,p_v
\,T^v_{sq}=\\
=\sum^n_{k=1}\sum^n_{q=1}\sum^n_{u=1}\sum^n_{m=1}\sum^n_{v=1}
\sum^n_{a=1}p_u\,T^u_{rk}\,P^k_a\,\tilde\nabla^aW^m\,
\frac{W^q\,p_m}{\Omega}\,p_v\,T^v_{sq}\,+\\
+\,\sum^n_{k=1}\sum^n_{q=1}\sum^n_{u=1}\sum^n_{v=1}\frac{p_u\,
T^u_{rk}\,W^k}{\Omega}\,\left(\,\shave{\sum^n_{a=1}\sum^n_{m=1}}
p_a\,\tilde\nabla^aW^m\,p_m\right)
\,\frac{p_v\,T^v_{sq}\,W^q}{\Omega}.
\endgathered
\tag13.22
$$
Last term in \thetag{13.22} is symmetric with respect to indices
$r$ and $s$. Therefore we can omit it adding to those denoted by
dots in formula for $C_{rs}$:
$$
\pagebreak
\hskip -2em
\gathered
C_{rs}\to C_{rs}+\ldots+\sum^n_{k=1}\sum^n_{q=1}\sum^n_{u=1}p_u\,T^u_{rk}
\,P^k_q\,B^q_s\,+\\
+\sum^n_{k=1}\sum^n_{q=1}\sum^n_{u=1}\sum^n_{m=1}\sum^n_{v=1}
\sum^n_{a=1}p_u\,T^u_{rk}\,P^k_q\,\tilde\nabla^qW^m\,
P^a_m\,p_v\,T^v_{sa}.
\endgathered
\tag13.23
$$
Now we can find transformation rule for left hand side of the equation
\thetag{11.8}. Using formula \thetag{13.23} and taking into account
equations \thetag{11.13} and \thetag{11.7}, we derive
$$
\hskip -2em
\sum^n_{r=1}\sum^n_{s=1}(C_{rs}-C_{sr})\,P^r_i\,P^s_j
\to\sum^n_{r=1}\sum^n_{s=1}(C_{rs}-C_{sr})\,P^r_i\,P^s_j.
\tag13.24
$$
Now, summarizing formulas \thetag{13.17}, \thetag{13.20}, \thetag{13.21},
and \thetag{13.24}, we see that the following theorem is proved.
\proclaim{Theorem 13.2} Additional normality equations \thetag{11.7},
\thetag{11.8}, \thetag{11.13} for Newtonian dynamical system \thetag{7.2}
are invariant under gauge transformations given by formulas \thetag{13.1},
\thetag{13.2}.
\endproclaim
    Theorems~13.1 and 13.2 mean that we can apply all normality equations
\thetag{8.21}, \thetag{11.7}, \thetag{11.8}, \thetag{11.13} to Newtonian
dynamical system written in connection-free form \thetag{2.3}. For this
purpose one should set $Q_i=\Theta_i$ in them, and one should choose
identically zero connection components $\Gamma^k_{ij}=0$, thus replacing
covariant derivatives $\nabla$ and $\tilde\nabla$ by corresponding partial
derivatives:
$$
\xalignat 2
&\nabla_{\!i}\to\frac{\partial}{\partial x^i},
&&\tilde\nabla^i\to\frac{\partial}{\partial p_i}.
\endxalignat
$$
However, in this form normality equations are not obviously coordinate
covariant. Each term in them loose transparent tensorial interpretation
provided by covariant derivatives. Generally speaking, we have a problem
of constructing present theory in coordinate-free and connection-free
form. This is separate problem, it will be studied in future papers.
\head
14. Basic example.
\endhead
   Let $H=H(x^1,\ldots,x^n,\,p_1,\ldots,p_n)$ be Hamilton function
for some Hamiltonian dynamical system. It is given by the following
well-known differential equations:
$$
\xalignat 2
&\hskip -2em
\dot x^i=\frac{\partial H}{\partial p_i},
&&\dot p_i=-\frac{\partial H}{\partial x^i}.
\tag14.1
\endxalignat
$$
Using some symmetric extended affine connection $\Gamma$, we can write
Hamilton equations \thetag{14.1} in terms of covariant derivatives
$\nabla$ and $\tilde\nabla$:
$$
\xalignat 2
&\hskip -2em
\dot x^i=\tilde\nabla^iH,
&&\nabla_{\!t}p_i=-\nabla_{\!i}H.
\tag14.2
\endxalignat
$$
Using the same Hamilton function $H=H(x^1,\ldots,x^n,\,p_1,\ldots,p_n)$
as in \thetag{14.2}, we define so called {\bf modified Hamiltonian
dynamical system}:
$$
\xalignat 2
&\hskip -2em
\dot x^i=\frac{\tilde\nabla^iH}{\dsize\sum^n_{s=1}p_s\,\tilde\nabla^sH},
&&\nabla_{\!t}p_i=-\frac{\nabla_{\!i}H}{\dsize\sum^n_{s=1}p_s\,\tilde
\nabla^sH}.
\tag14.3
\endxalignat
$$
Surely we should choose $H$ such that denominator in \thetag{14.3} is
non-zero. Dynamical system \thetag{14.3} is an example (is special
case) of Newtonian dynamical system \thetag{7.2} with vector field
$\bold V$ and covector field $\bold Q$ given by formulas
$$
\xalignat 2
&\hskip -2em
V^i=\frac{\tilde\nabla^iH}{\dsize\sum^n_{s=1}p_s\,\tilde\nabla^sH},
&&Q_i=-\frac{\nabla_{\!i}H}{\dsize\sum^n_{s=1}p_s\,\tilde
\nabla^sH}.
\tag14.4
\endxalignat
$$
Substituting \thetag{14.4} into \thetag{6.1} and \thetag{6.2}, we find
$$
\xalignat 2
&\hskip -3em
W^i=-V^i,
&&\Omega=-1.
\tag14.5
\endxalignat
$$
\proclaim{Theorem 14.1} Modified Hamiltonian dynamical system
\thetag{14.3} is an example of Newtonian dynamical system
\thetag{7.2} admitting normal shift of hypersurfaces in the sense
of definition~1.5.
\endproclaim
We can write modified dynamical system \thetag{14.3} in connection-free
form:
$$
\xalignat 2
&\hskip -2em
\dot x^i=\frac{\dfrac{\partial H}{\partial p_i}
\vphantom{\vrule height 1pt depth 12pt}}
{\dsize\sum^n_{s=1}p_s\,\dfrac{\partial H}{\partial p_s}},
&&\dot p_i=-\frac{\dfrac{\partial H}{\partial x^i}
\vphantom{\vrule height 1pt depth 10pt}}
{\dsize\sum^n_{s=1}p_s\,\dfrac{\partial H}{\partial p_s}}.
\tag14.6
\endxalignat
$$
\proclaim{Theorem 14.2} Modified Hamiltonian dynamical system
\thetag{14.6} is an example of Newtonian dynamical system
\thetag{2.3} admitting normal shift of hypersurfaces in the sense
of definition~1.5.
\endproclaim
Differential equations \thetag{14.3} and \thetag{14.6} are different
representations of the same dynamical system. Therefore theorems~14.1
and 14.2 formulate the same result. This result is not new. It was
obtained in paper \cite{19}. In present paper it is built into framework
of more general theory and forms basic example for this theory.
It proves that class of dynamical systems considered in this theory is
not empty, and, moreover, this class is sufficiently large.
\head
15. Dedication.
\endhead
   This paper is dedicated to my mother F. M. Sharipova, who taught
Mathematics for many years in School no\.~18 of Karakul, Bukhara
region, Uzbekistan. She was best in recognizing various constellations
on the sky and knew great many of them. In dark south nights in Summer,
when I was 10 or even younger, she often told me about stars, comets,
planets, and other thing. Possibly this was why later on I have chosen
Natural Sciences and Mathematics for my profession. Bright image of
my mother is ever kept in my memory.\par


\newpage
\Refs
\ref\no 1\by Boldin~A.~Yu\., Sharipov~R.~A.\book Dynamical systems
accepting the normal shift\publ Preprint No\.~0001-M of Bashkir State
University\publaddr Ufa\yr April, 1993
\endref
\ref\no 2\by Boldin~A.~Yu\., Sharipov~R.~A.\paper Dynamical systems
that admit normal shift\jour Teoret\. Mat\. Fiz\. \vol 97\issue 3
\yr 1993\pages 386--395\moreref see also chao-dyn/9403003 in Electronic
Archive at LANL\footnotemark
\endref
\footnotetext{Electronic Archive at Los Alamos National Laboratory of
USA (LANL). Archive is accessible through Internet 
{\bf http:/\negskp/arXiv\.org}, it has mirror site 
{\bf http:/\negskp/ru\.arXiv\.org} at the Institute for Theoretical
and Experimental Physics (ITEP, Moscow).}
\adjustfootnotemark{-1}
\ref\no 3\by Boldin~A.~Yu\., Sharipov~R.~A.\paper Dynamical systems
that admit normal shift\jour Dokl\. Akad\. Nauk \vol 334\yr 1994
\issue 2\pages 165--167
\endref
\ref\no 4\by Boldin~A.~Yu\., Sharipov~R.~A.\paper Multidimensional
dynamical systems that admit normal shift\jour Teoret\. Mat\. Fiz\. 
\vol 100\issue 2\yr 1994\pages 264--269\moreref see also
patt-sol/9404001 in Electronic Archive at LANL
\endref
\ref\no 5\by Boldin~A.~Yu\., Dmitrieva~V.~V., Safin~S.~S., Sharipov~R.~A.
\paper Dynamical systems on Riemannian manifolds that admit normal shift
\jour Teoret\. Mat\. Fiz\. \yr 1995\vol 105\issue 2 \pages 256--266
\moreref\inbook see also in book ``{Dynamical systems accepting
the normal shift}''\publ Bashkir State University\publaddr Ufa\yr 1994
\pages 4--19\moreref and see hep-th/9405021 in Electronic Archive at
LANL
\endref
\ref\no 6\by Boldin~A.~Yu\., Bronnikov~A.~A., Dmitrieva~V.~V.,
Sharipov~R.~A.\paper Complete normality conditions for dynamical
systems on Riemannian manifolds\jour Teoret\. Mat\. Fiz\. \yr 1995
\vol 103\issue 2\pages 267--275\moreref\inbook see also in book
``{Dynamical systems accepting the normal shift}''\publ Bashkir
State University\publaddr Ufa\yr 1994\pages 20--30\moreref and
see astro-ph/9405049 in Electronic Archive at LANL
\endref
\ref\no 7\by Sharipov~R.~A.\book Dynamical systems admitting the normal
shift\publ Thesis for the degree of Doctor of Sciences in Russia\publaddr
Ufa\yr 1999\moreref English version of thesis is submitted to Electronic
Archive at LANL, see archive file math.DG/0002202
\endref
\ref\no 8\by Boldin~A.~Yu\.\book Two-dimensional dynamical systems
admitting the normal shift\publ Thesis for the degree of Candidate of
Sciences in Russia\yr 2000\moreref English version of thesis is
submitted to Electronic Archive at LANL, see archive file math.DG/0011134
\endref
\ref\no 9\by Boldin~A\.~Yu\.\paper On the self-similar solutions of 
normality equation in two-dimensional case\inbook in book ``{Dynamical
systems accepting the normal shift}''\publ Bashkir State University
\publaddr Ufa\yr 1994\pages 31--39\moreref see also patt-sol/9407002
in Electronic Archive at LANL
\endref
\ref\no 10\by Boldin~A.~Yu\., Sharipov~R.~A.\paper On the solution
of normality equations in the dimension\linebreak $n\geqslant 3$\jour
Algebra i Analiz\vol 10\yr 1998\issue 4\pages 37--62\moreref see also
solv-int/9610006 in Electronic Archive at LANL
\endref
\ref\no 11\by Sharipov~R.~A.\paper Newtonian normal shift in
multidimensional Riemannian geometry\jour Mat\. Sb\.\vol 192
\issue 6\yr 2001\pages 105--144\moreref\jour see also 
math.DG/0006125 in Electronic Archive at LANL
\endref
\ref\no 12\by Sharipov~R.~A.\paper Newtonian dynamical systems
admitting normal blow-up of points\jour Zap\. sem\. POMI
\vol 280\yr 2001\pages 278--298\moreref see also proceeding
of Conference organized by R.~S.~Saks in Ufa, August 2000,
pp\.~215-223, and math.DG/0008081 in Electronic Archive at LANL
\endref
\ref\no 13\by Sharipov~R.~A.\paper On the solutions of weak
normality equations in multidimensional case\jour Paper
math.DG/0012110 in Electronic Archive at LANL\yr 2000
\endref
\ref\no 14\by Sharipov~R.~A. Global geometric structures associated
with dynamical systems admitting normal shift of hypersurfaces in
Riemannian manifolds\jour International Journ\. of Mathematics and
Math\. Sciences \vol 30\issue 9\yr 2002\pages 541--558\moreref\paper
{\rm see also} First problem of globalization in the theory of dynamical
systems admitting the normal shift of hypersurfaces\jour math.DG/0101150
in Electronic Archive at LANL\yr 2001
\endref
\ref\no 15\by Sharipov~R.~A.\paper Second problem of globalization
in the theory of dynamical systems admitting the normal shift of
hypersurfaces\jour Paper math.DG/0102141 in Electronic Archive at
LANL\yr 2001
\endref
\ref\no 16\by Sharipov~R.~A.\paper A note on Newtonian, Lagrangian,
and Hamiltonian dynamical systems in Riemannian manifolds\jour Paper
math.DG/0107212 in Electronic Archive at LANL\yr 2001
\endref
\ref\no 17\by Sharipov~R.~A.\paper Dynamic systems admitting the normal
shift and wave equations\jour Teoret\.~Mat\. Fiz\. \vol 131\issue 2
\pages 244--260\yr 2002\moreref see also math.DG/0108158 in Electronic
Archive at LANL
\endref
\ref\no 18\by Fedoryuk~M.~V. \paper The equations with fast oscillating
solutions\inbook Summaries of Science and Technology. Modern problems
of Mathematics. Fundamental Researches. Vol. 34\yr 1988\publ VINITI
\publaddr Moscow
\endref
\ref\no 19\by Sharipov~R.~A.\paper Normal shift in general
Lagrangian dynamics\jour Paper math.DG/0112089 in Electronic
Archive at LANL\yr 2001
\endref
\ref\no 20\by Sharipov~R.~A.\paper Comparative analysis for pair of
dynamical systems, one of which is Lagrangian\jour Paper math.DG/0204161
in Electronic Archive at LANL\yr 2002
\endref
\ref\no 21\by Arnold~V.~I.\book Mathematical methods of classical
mechanics\publ Nauka publishers\publaddr Moscow\yr 1979
\endref
\ref\no 22\by Filippov~V.~M., Savchin~V.~M., Shorohov~S.~G.
\paper Variational principles for non-potential operators
\inbook Modern problems in mathematics. Recent achievements
\vol 40\yr 1992\pages 3--176\publ VINITI\publaddr Moscow
\endref
\endRefs
\enddocument
\end